\documentclass[11pt,a4paper]{amsart}
\usepackage[T1]{fontenc}
\usepackage[utf8]{inputenc}
\usepackage{lmodern}
\usepackage[a4paper,margin=29mm]{geometry}
\usepackage{amsmath,amssymb,mathtools,mathrsfs}
\usepackage{enumitem}
\usepackage{needspace}
\usepackage[expansion=false]{microtype}
\usepackage{hyperref}
\hypersetup{colorlinks=true,linkcolor=blue,citecolor=blue,urlcolor=blue,
 pdftitle={Hessian sectional curvatures of conjugate symmetric statistical manifolds},
 pdfauthor={Hikozo Kobayashi}}
\numberwithin{equation}{section}
\newtheorem{Thm}{Theorem}[section]
\newtheorem{Prop}[Thm]{Proposition}
\newtheorem{Lem}[Thm]{Lemma}
\newtheorem{Cor}[Thm]{Corollary}
\theoremstyle{definition}
\newtheorem{Def}[Thm]{Definition}
\newtheorem{Rem}[Thm]{Remark}
\newtheorem{Ex}[Thm]{Example}
\newtheorem{Question}[Thm]{Question}
\newcommand{\RR}{\mathbb R}

\newcommand{\calQ}{\mathcal Q}

\newcommand{\Aut}{\mathrm{Aut}}
\newcommand{\MLStat}{\mathcal{ML}\mathrm{Stat}}

\newcommand{\Hess}{\mathrm{Hess}}
\newcommand{\LR}{\langle, \rangle}
\newcommand{\gee}{\mathfrak{g}}
\newcommand{\LStat}{\mathcal{L}\mathrm{Stat}}
\newcommand{\calK}{\mathcal K}
\newcommand{\STAT}{\mathop{\mathbf{Stat}}}
\newcommand{\HS}{\mathop{\mathbf{HS}}}
\newcommand{\Kflat}{K^\flat}

\DeclareMathOperator{\PHC}{PHC}
\DeclareMathOperator{\tr}{tr}
\DeclareMathOperator{\pr}{pr}
\DeclareMathOperator{\id}{id}

\setlist[enumerate]{itemsep=3pt,topsep=5pt}
\allowdisplaybreaks[1]
\title[Hessian sectional curvatures of conjugate symmetric statistical manifolds]{Hessian sectional curvatures of conjugate symmetric statistical manifolds}
\author{Hikozo Kobayashi}
\thanks{The author was supported by JST SPRING, Grant Number JPMJSP2132.}
\address{Mathematics Program, Graduate School of Advanced Science and
 Engineering, Hiroshima University, 1-3-1 Kagamiyama, Higashi-Hiroshima
 City, Hiroshima 739-8526, Japan}
\email{hikozo-kobayashi@hiroshima-u.ac.jp}
\begin{document}
\begin{abstract}
Conjugate symmetric statistical manifolds generalize Hessian manifolds.
Shima introduced Hessian sectional curvature in 1995 as an invariant for studying Hessian structures.
We extend this invariant to conjugate symmetric statistical manifolds and study related curvature conditions and invariants.
Our main results are as follows.
For connected Hessian manifolds of dimension at least two, we characterize constant Hessian sectional curvature (CHC) by pointwise orthogonal invariance.
We call conjugate symmetric statistical manifolds satisfying
this invariance condition weakly CHC.
We classify connected weakly CHC statistical manifolds of dimension at least two whose cubic forms are parallel with respect to the Levi-Civita connection.
This classification characterizes the non-Hessian statistical
structures induced on the classical hyperbolic affine spheres
$xyz=c>0$ in the positive octant.
We also introduce principal Hessian curvatures and show that they form a complete invariant for left-invariant Hessian structures on $\mathbb{R}^n$ up to statistical isomorphism.
\end{abstract}
\maketitle
\begingroup
\footnotesize
\tableofcontents
\endgroup

\section{Introduction}\label{sec:introduction}

The construction and study of invariants play an important role in differential geometry.
In this paper, we focus on Hessian manifolds and the broader class of conjugate symmetric statistical manifolds, and study their invariants.

We first consider Hessian sectional curvature on Hessian manifolds and conditions for its constancy.
A \emph{Hessian manifold} is a manifold equipped with a Riemannian metric $g$ and a flat torsion-free affine connection $\nabla$ such that $\nabla g$ is symmetric.
Hessian manifolds were introduced and systematically studied by Shima~\cite{Shima_1976, Shima_1980, Shima_1995, Shima} and are closely related to information geometry (see, for example, \cite{Amari-Nagaoka_2000, Shima}).

Shima~\cite{Shima_1995} introduced Hessian sectional curvature as an invariant of Hessian structures.
We briefly recall its definition (see Sections~\ref{subsec:Hessian_curvature_tensor} and \ref{subsec:Hessian_sectional_principal} for details).
Let $(M,g,\nabla)$ be a Hessian manifold.
Then, a symmetric bilinear form $\calQ$ is defined on $S^2TM$, the bundle of contravariant symmetric tensors of degree $2$ on $M$.
For $p\in M$ and $0\ne\nu\in S^2T_pM$, the \emph{Hessian sectional curvature in the direction $\nu$} is defined by
\[
q_p(\nu):=\frac{\calQ_p(\nu,\nu)}{ \langle \nu,\nu \rangle_{g_p} },
\]
where $\LR_{g_p}$ is the inner product on $S^2T_pM$ induced by $g_p$.

If there exists a constant $c \in \RR$ such that $q_p(\nu)=c$ for every $p\in M$ and $0\ne\nu\in S^2T_pM$, then $(M,g,\nabla)$ is said to be \emph{of constant Hessian sectional curvature $c$}~\cite{Shima_1995}, or \emph{of CHC $c$} for short~\cite{FK_2013}.
The classification and submanifold theory of Hessian manifolds of CHC have attracted interest and been studied (see, for example,
\cite{Chen-Mihai-Mihai_2019, Furuhata_2009, Furuhata_2011,
Furuhata-Hasegawa-Satoh_2022, FK_2013, Mihai-Mihai_2020,
Mihai_2020, Sakamoto_2026, Shima_1995, Siddiqui-Ahmad-Ozel_2019}).

If a Hessian manifold is of CHC, then its Hessian sectional curvature at each point is invariant under the action of the orthogonal group of the tangent space. 
In this paper, we show that the converse holds for connected Hessian manifolds of dimension at least $2$. 
This is our first main result:

\begin{Thm}[cf.~Theorem~\ref{thm:hessian_weak}]\label{thm:main-1_intro}
Let $(M^n,g,\nabla)$ be a connected Hessian manifold of dimension $n\geq2$.
Then, the following conditions are equivalent.
\begin{enumerate}[label=\textup{(\roman*)}]
    \item\label{item:main-1_intro_item1}
    There exists $c \in \RR$ such that $(M,g,\nabla)$ is of CHC $c$.
    \item\label{item:main-1_intro_item2}
   For every $p\in M$ and every $0\ne\nu\in S^2T_pM$,
    \[
    q_p(u\cdot\nu)=q_p(\nu)
    \qquad \bigl(u\in O(T_pM,g_p)\bigr).
    \]
\end{enumerate}
\end{Thm}

On the other hand, in dimension $1$, this equivalence does not hold in general, even when $M$ is connected (cf.~Remark~\ref{rem:dimension_one}).

Next, we extend the notion of Hessian sectional curvature to conjugate symmetric statistical manifolds (see Definitions~\ref{def:statistical} and~\ref{def:conjugate_symmetric} for the definition).
This class has also been studied from various viewpoints, including information geometry, the theory of statistical manifolds, and affine hypersurface theory~\cite{FIK, Inoguchi-Ohno_2024, IKO_2026, KOOT_2025, Kobayashi-Okuda_2025, Kobayashi-Ohno_2025, Kurose-1990, Kurose_1994,
L-SM, Matsuzoe-Takeuchi-Amari_2006, Min-Ri-Kwak-An_2015, Noguchi_1992,Nomizu-Sasaki,Opozda_2015,Opozda_2016, Opozda_2019, Opozda_2021, 
Opozda_2022,Opozda_2025, Takeuchi-Amari_2005}.
It is therefore important to introduce invariants that are useful for the classification of manifolds in this class and for the study of their geometric properties.

We asked whether the notion of Hessian sectional curvature could be
extended to conjugate symmetric statistical manifolds.
Our second main result is the following theorem,
which gives an affirmative answer to this question.

\begin{Thm}[see~Section~\ref{subsec:Hessian_sectional_principal} for details]
Hessian sectional curvature can also be defined for conjugate symmetric statistical manifolds and is invariant under statistical isomorphisms (see Definition~\ref{def:statistical_isomorphism}).
For Hessian manifolds, it agrees with that introduced by Shima~\cite{Shima_1995}.
\end{Thm}

This extension also allows us to introduce the notion of CHC for conjugate symmetric statistical manifolds.
We study sufficient conditions for CHC (see Sections~\ref{subsec:K_curvature_Schur} and~\ref{subsec:Hessian_weak_CHC}) and construct non-Hessian examples (see Proposition~\ref{prop:warped_CHC}).

Moreover, motivated by condition~\ref{item:main-1_intro_item2} of Theorem~\ref{thm:main-1_intro}, we introduce the following notion.

\begin{Def}[cf.~Definition~\ref{def:weak_CHC}]
A conjugate symmetric statistical manifold $(M,g,\nabla)$ is said to be  \emph{weakly CHC} if
\[
    q_p(u\cdot\nu)=q_p(\nu)
\]
for every $p\in M$, $0\ne\nu\in S^2T_pM$, and $u\in O(T_pM,g_p)$.
\end{Def}

Theorem~\ref{thm:main-1_intro} states that weakly CHC and CHC are equivalent for connected Hessian manifolds of dimension at least $2$.
This naturally leads to the following question.

\begin{Question}
Are the weakly CHC and CHC conditions also equivalent for connected conjugate symmetric statistical manifolds of dimension at least $2$?
\end{Question}

The answer to this question is negative in general.
A counterexample is provided by the following classical affine surface.
For $c>0$, let
\[
    \Sigma_c:=\{(x,y,z)\in(\RR_{>0})^3\mid xyz=c\}.
\]
Let $g_c$ and $\nabla_c$ denote its affine metric and Blaschke connection, respectively.
Here, the direction of the affine normal field is chosen so that $g_c$ is positive definite.
Then $(\Sigma_c,g_c,\nabla_c)$ is a non-Hessian conjugate symmetric
statistical manifold that is weakly CHC but not of CHC (see Section~\ref{subsec:WCHC_Classification} for details).

The surface $\Sigma_c$ is a classical and important example
of a hyperbolic affine sphere (cf.~\cite{Calabi_1972, Dumas-Wolf_2015, Li-Penn_1988, Loftin_2010, Loftin-McIntosh_2016, Magid-Nomizu_1989, Magid-Ryan_1990, Tzitzeica_1908, Vrancken-Li-Simon_1991}).
In this paper, we classify weakly CHC conjugate symmetric
statistical manifolds whose cubic forms are parallel with
respect to the Levi-Civita connection and give a local
characterization of the statistical structures associated
with these surfaces, as follows.
This is our third main result:

\begin{Thm}[cf.~Theorem~\ref{thm:parallel_classification}]
\label{thm:main-2_intro}
Let $(M^n,g,\nabla)$ be a connected conjugate symmetric statistical manifold of dimension $n\geq2$, 
and assume that its cubic form $C:=\nabla g$ is parallel with respect to the Levi-Civita connection $\nabla^g$.
Then $(M,g,\nabla)$ is weakly CHC if and only if one of the following holds:
\begin{enumerate}[label=\textup{(\roman*)}]
    \item $\nabla=\nabla^g$.
    \item $n=2$, and there exists $c>0$ such that $(M,g,\nabla)$ is
locally statistically isomorphic to $(\Sigma_c,g_c,\nabla_c)$
in the sense of Definition~\ref{def:locally_statistical_isomorphism}.
\end{enumerate}
In the latter case, $c$ is uniquely determined, and if $M$ is
simply connected and $g$ is complete, then $(M,g,\nabla)$ is
globally statistically isomorphic to $(\Sigma_c,g_c,\nabla_c)$.
\end{Thm}

Finally, we introduce the eigenvalues of the self-adjoint operator associated with Hessian sectional curvature as \emph{principal Hessian curvatures}.
We then investigate to what extent they distinguish statistical structures.

Vaisman defined \emph{principal conical curvatures} in the Hessian setting
in~\cite[Open Problem in Appendix]{Vaisman_2014}.
The following definition extends this notion to conjugate symmetric
statistical manifolds under the name principal Hessian curvatures.

\begin{Def}[cf.~Definition~\ref{def:Hessian_curvatures}]
Let $(M,g,\nabla)$ be a conjugate symmetric statistical manifold and $p\in M$.
We say that each eigenvalue of the self-adjoint operator
$\hat{Q}_p:S^2T_pM\to S^2T_pM$ defined by
\[
    \calQ_p(\nu,\mu)
    =\langle \hat{Q}_p(\nu), \mu \rangle_{g_p}
    \qquad \bigl(\nu,\mu\in S^2T_pM\bigr)
\]
is a \emph{principal Hessian curvature at $p$}.
\end{Def}

The principal Hessian curvatures retain only part of the information contained in Hessian sectional curvature.
Together with the corresponding eigenspaces, they determine Hessian sectional curvature (see Section~\ref{subsec:Hessian_sectional_principal} for details).
As with Hessian sectional curvature, the principal Hessian curvatures, counted with multiplicity, are preserved under statistical isomorphisms (Theorem~\ref{thm:HSC-def-on-CS}).

As a concrete application of this invariant, we consider left-invariant Hessian structures on the abelian Lie group $\RR^n$.
For every left-invariant Hessian structure $(g,\nabla)$ on $\RR^n$, both $g$ and $\nabla$ are flat, but such structures need not be statistically isomorphic to one another.
We explicitly compute the principal Hessian curvatures of all left-invariant Hessian structures and obtain the following result.
This is the fourth main result of this paper:

\begin{Thm}[cf.~Theorem~\ref{thm:euclidean_spectrum}]\label{thm:main-3_intro}
The principal Hessian curvatures, counted with multiplicity, form a complete invariant for left-invariant Hessian structures on the abelian Lie group $\RR^n$ up to statistical isomorphism.
\end{Thm}

We also present two conjugate symmetric statistical manifolds that are not Hessian and can be distinguished by their principal Hessian curvatures, although they have the same values of several other invariants (see~Example~\ref{ex:PHC_distinguishes} for details).
These results suggest that principal Hessian curvatures may be useful for the classification of statistical structures and the study of their moduli spaces.

For Hessian manifolds, Hessian sectional curvature is closely related to the holomorphic sectional curvature of the natural K\"ahler metric on the tangent bundle (see~\cite{Shima} for details).
For related studies, see also~\cite{Khan-Zhang_2020,Khan-Zhang-Zheng_2025,Khan-Zheng_2026, Figueiredo-Molitor_2026, Satoh_2007,Zhang-Khan_2020}.
The relationship between the curvature conditions considered in this paper and the almost K\"ahler structures induced on the tangent bundles of more general statistical manifolds will be studied in detail in our subsequent paper.

This paper is organized as follows.
After reviewing the preliminaries in Section~\ref{sec:preliminaries}, we introduce Hessian sectional curvature and principal Hessian curvatures for conjugate symmetric statistical manifolds in Section~\ref{sec:curvature}.
In Section~\ref{sec:conditions}, we study the CHC and weakly CHC conditions, give explicit examples, and prove Theorem~\ref{thm:main-1_intro}.
In Section~\ref{sec:euclidean}, we consider left-invariant Hessian structures on $\RR^n$ and prove Theorem~\ref{thm:main-3_intro}.
In Section~\ref{sec:parallel}, we state the classification of weakly CHC conjugate symmetric statistical manifolds with parallel cubic forms and the characterization corresponding to Theorem~\ref{thm:main-2_intro}, with proofs given in Appendix~\ref{app:parallel}.

\section{Preliminaries}\label{sec:preliminaries}

Throughout this paper, all manifolds and tensor fields are smooth, and
all metrics are positive definite. We write $M^n$ for an $n$-dimensional
manifold $M$. 
For a vector bundle $E$, we denote by $\Gamma(E)$ the space
of its smooth sections and by $S^kE$ its $k$th symmetric tensor
power. We also use the notation $S^kE$ for vector spaces.
For $\eta_1,\ldots,\eta_k\in\Gamma(E)$, we use the convention
\[
 \eta_1\odot\cdots\odot\eta_k
 =\frac{1}{k!}\sum_{\sigma\in\mathfrak{S}_k}
 \eta_{\sigma(1)}\otimes\cdots\otimes\eta_{\sigma(k)},
\]
where $\mathfrak{S}_k$ denotes the symmetric group on $k$ letters. We omit the symbol $\odot$ when no confusion can arise. The curvature tensor of an affine connection $\nabla$ is denoted by $R^\nabla$, with the convention
\[
 R^\nabla(X,Y)Z
 =\nabla_X\nabla_YZ-\nabla_Y\nabla_XZ-\nabla_{[X,Y]}Z.
\]

For a Riemannian metric $g$ on $M$, we denote its Levi-Civita connection
by $\nabla^g$ and set $R^g:=R^{\nabla^g}$.
On tensor bundles and their subbundles, we write $\LR_g$ for the
fiber metric induced by $g$ and $\LR_{g_p}$ for the corresponding
inner product on the fiber over $p\in M$.
For a tensor field $T$, we denote the resulting norm of $T_p$ by
$|T_p|$ and the function $p\mapsto |T_p|$ by $\lVert T\rVert_g$.

In this section, we recall some basic notions and facts
concerning statistical manifolds that will be used in this paper.
For further details, see~\cite{Amari-Nagaoka_2000, L-SM, Noguchi_1992, Opozda_2015} for example.

\subsection{Statistical manifolds and their isomorphisms}\label{subsec:Stat-str}

We first recall statistical structures and the notion of statistical isomorphism.

\begin{Def}[cf.~\cite{Kurose-1990}]\label{def:statistical}
A \emph{statistical structure} on a manifold $M$ is a pair $(g,\nabla)$
consisting of a Riemannian metric and a torsion-free affine connection
such that $C:=\nabla g$ is a symmetric $(0,3)$-tensor field.
The triple $(M,g,\nabla)$ is called a \emph{statistical manifold},
$\nabla$ a \emph{statistical connection} on $(M,g)$, and
$C$ its \emph{cubic form}.
\end{Def}

For a statistical manifold $(M,g,\nabla)$, its \emph{difference tensor} $K:=\nabla-\nabla^g$ satisfies
\begin{equation}\label{eq:K-C}
 K(X,Y)=K(Y,X),\qquad
 \Kflat(X,Y,Z):=g(K(X,Y),Z)=-\frac12 C(X,Y,Z).
\end{equation}
In particular, $\Kflat$ is symmetric, and the endomorphism
$K_X:Y\mapsto K(X,Y)$ is self-adjoint with respect to $g$.
Conversely, a metric $g$ and a symmetric $(0,3)$-tensor field $C$
determine $K$ by~\eqref{eq:K-C} and hence a statistical connection
$\nabla=\nabla^g+K$.
We therefore also write $(g,C)$ for a statistical structure
and $(M,g,C)$ for a statistical manifold.
The structure is \emph{trivial} if $K=0$, equivalently, if
$\nabla=\nabla^g$.

The \emph{dual connection} $\nabla^*$ is defined by
\begin{equation}\label{eq:def_dual-conn}
    Xg(Y,Z)=g(\nabla_XY,Z)+g(Y,\nabla_X^*Z).
\end{equation}
It is torsion-free and satisfies
\begin{equation}\label{eq:dual}
 \nabla^*=\nabla^g-K,\qquad \nabla^*g=-C.
\end{equation}

We now introduce the following notions of isomorphism.

\begin{Def}[{cf.~\cite[Definition~1.2]{Satoh-2019}; see also \cite[Section~9]{L-SM}}]
\label{def:statistical_isomorphism}
A \emph{statistical isomorphism} from $(M,g,\nabla)$ to
$(M',g',\nabla')$ is a diffeomorphism $f:M\to M'$ satisfying
$f^*g'=g$ and $f^*\nabla'=\nabla$. 
We say that $(M,g,\nabla)$ and $(M',g',\nabla')$ are
\emph{statistically isomorphic} if there exists a statistical
isomorphism between them.
\end{Def}

A diffeomorphism is a statistical isomorphism if and only if it
preserves $g$ and $K$, or equivalently, $g$ and $C$.

\begin{Def}\label{def:STAT_equivalence}
For two statistical structures $(g,\nabla)$ and $(g',\nabla')$ on a
manifold $M$, we write
\[
 (g,\nabla)\sim_{\STAT}(g',\nabla')
\]
if there exists a statistical isomorphism from $(M,g,\nabla)$ to
$(M,g',\nabla')$.
\end{Def}

This defines an equivalence relation on the set of statistical
structures on $M$.

We conclude this subsection by specifying what we mean by two
statistical manifolds being locally statistically isomorphic.

\begin{Def}\label{def:locally_statistical_isomorphism}
We say that $(M,g,\nabla)$ and $(M',g',\nabla')$ are
\emph{locally statistically isomorphic} if every point of $M$
has an open neighborhood statistically isomorphic to an open
subset of $M'$, and the same holds with $M$ and $M'$ interchanged.
\end{Def}

\subsection{Conjugate symmetry and related curvature conditions}\label{subsec:class_structures}

We introduce the curvature conditions used in this paper
and recall some related results.

Let $(M, g, \nabla)$ be a statistical manifold. The following curvature identities are well known; see, for example, \cite[Section~3]{Opozda_2015}:
\begin{align}
 R^\nabla(X,Y)Z
 &=R^g(X,Y)Z+(\nabla^g_XK)(Y,Z)-(\nabla^g_YK)(X,Z)
      +[K_X,K_Y]Z,\label{eq:curvature_decomposition}\\
 g(R^\nabla(X,Y)Z,W)
 &=-g(Z,R^{\nabla^*}(X,Y)W).\label{eq:dual_curvature}
\end{align}
Here $[K_X,K_Y]=K_XK_Y-K_YK_X$.

\begin{Def}[cf.~\cite{L-SM,Shima}]\label{def:conjugate_symmetric}
A statistical manifold is \emph{conjugate symmetric} if
$R^\nabla=R^{\nabla^*}$ (see~\eqref{eq:def_dual-conn} for the definition of $\nabla^*$). It is \emph{Hessian} if $R^\nabla=0$.
\end{Def}

Conjugate symmetry is equivalent to
\begin{equation}\label{eq:CS_Codazzi}
 (\nabla^g_XK)(Y,Z)=(\nabla^g_YK)(X,Z),
\end{equation}
or to total symmetry of $\nabla^gC$, equivalently, of
$\nabla^g\Kflat$
(see~\cite{Noguchi_1992}).
In particular, $\nabla^g C=0$ is a sufficient condition for conjugate symmetry, and every trivial statistical structure satisfies this condition since its cubic form vanishes.
By~\eqref{eq:dual_curvature}, conjugate symmetry is also
equivalent to $g$-skew-adjointness of $R^\nabla(X,Y)$.
For conjugate symmetric statistical structures,
\eqref{eq:curvature_decomposition} reduces to
\begin{equation}\label{eq:conjugate_curvature}
 R^\nabla(X,Y)Z=R^g(X,Y)Z+[K_X,K_Y]Z.
\end{equation}

We recall the following consequence of a result of
Noguchi~\cite[Proposition~5.1]{Noguchi_1992}
for conjugate symmetric statistical manifolds.

\begin{Prop}\label{prop:pure_trace_Codazzi}
Let $(M^n,g,\nabla)$ be a conjugate symmetric statistical
manifold. Suppose that
\[
 \Kflat(X,Y,Z)
 =\theta(X)g(Y,Z)+\theta(Y)g(X,Z)+\theta(Z)g(X,Y)
\]
for a $1$-form $\theta$ on $M$ and all vector fields $X,Y,Z$.
Then there exists a function $\sigma$ on $M$ such that
$\nabla^g\theta=\sigma g$.
\end{Prop}

Put
\[
 R_0(X,Y)Z:=g(Y,Z)X-g(X,Z)Y.
\]
A statistical structure has \emph{constant curvature $k$} if
$R^\nabla=kR_0$ for a constant $k \in \RR$ (\cite{Kurose-1990}).
Thus Hessian structures are precisely statistical structures of constant
curvature $0$.
Every statistical structure of constant curvature is conjugate
symmetric by~\eqref{eq:dual_curvature}.
For examples and further background on these classes, see,
for example,~\cite{Amari-Nagaoka_2000,KO-CC,Kurose-1990,L-SM,
Noguchi_1992,Opozda_2015,Shima}.

We denote by $[K,K]$ the $(1,3)$-tensor defined by
\[
    [K,K](X,Y)Z := [K_X,K_Y]Z.  
\]
Following~\cite{Opozda_2016}, we say that a statistical structure
(or a statistical manifold) has
\emph{constant sectional $\mathsf K$-curvature $c$ at $p$},
where $p\in M$ and $c\in\RR$, if
\[
 [K,K]_p=c(R_0)_p.
\]
We say that it has
\emph{pointwise constant sectional $\mathsf K$-curvature $\rho$} if
\[
 [K,K]=\rho R_0
\]
for some smooth function $\rho$ on $M$.
When $\rho$ is constant, we omit ``pointwise.''

We recall below some results of Opozda~\cite{Opozda_2016}
on constant sectional $\mathsf K$-curvature.

\begin{Prop}[{\cite[Proposition~3.3]{Opozda_2016}}]
\label{prop:Opozda_K_Schur}
Let $(M^n,g,\nabla)$ be a connected conjugate symmetric
statistical manifold with $n\geq3$.
If the structure has pointwise constant sectional
$\mathsf K$-curvature $\rho$, then $\rho$ is constant on $M$.
\end{Prop}

\begin{Prop}[{cf.~\cite[Lemma~3.21]{Opozda_2016}}]
\label{prop:Opozda_negative_K}
Let $(M^n,g,\nabla)$ be a statistical manifold with
constant sectional $\mathsf K$-curvature $c<0$
at a point $p\in M$.
If an endomorphism $J:T_pM\to T_pM$ is skew-adjoint
with respect to $g_p$ and satisfies
\[
 J(K_p(X,Y))-K_p(JX,Y)-K_p(X,JY)=0
 \qquad(X,Y\in T_pM),
\]
then $J=0$.
\end{Prop}

\begin{Rem}
Proposition~\ref{prop:Opozda_negative_K} remains valid if
the sectional $\mathsf K$-curvature at $p$ is negative
but not necessarily constant; see~\cite[Lemma~3.21]{Opozda_2016}.
\end{Rem}

For further details on sectional $\mathsf K$-curvature,
see~\cite{Opozda_2016}.

For a local $g$-orthonormal frame $\{e_1,\ldots,e_n\}$, set
\begin{equation}\label{eq:trace-g_K}
 \tr_gK:=\sum_{i=1}^n K(e_i,e_i).
\end{equation}
This vector field is independent of the choice of frame. The statistical structure is called \emph{trace-free} if $\tr_gK=0$.
By~\cite[Section~3]{Opozda_2015}, we have
\begin{equation}\label{eq:trace_free_equivalences}
 \tr K_X=g(\tr_gK,X),\qquad
 \tr_g C(X,\cdot,\cdot)=-2g(\tr_gK,X).
\end{equation}
Thus the structure is trace-free if and only if $\tr K_X=0$ for every vector field $X$.

The vector field $\tr_g K$ also appears in the following scalar curvature identity.
Let $\operatorname{Scal}^{\nabla}:=\tr_g\operatorname{Ric}^{\nabla}$,
and let $\operatorname{Scal}^{g}$ denote the scalar curvature of $g$.
By~\cite[Equation~(19)]{Opozda_2022}, we have
\begin{equation}\label{eq:stat_scalar_norm}
 \operatorname{Scal}^{\nabla}-\operatorname{Scal}^{g}
 =\lVert\tr_gK\rVert_g^2-\lVert K\rVert_g^2.
\end{equation}

We will also use the standard Ricci identity
(cf.~\cite[Chapter~III, Proposition~2.12 and Theorem~5.1]{Kobayashi-Nomizu_I})
for a torsion-free affine connection $\nabla$ and a
$(1,2)$-tensor field $K$:
\begin{equation}\label{eq:Ricci_commutation}
 (\nabla_U\nabla K)(X,Y,Z)-(\nabla_X\nabla K)(U,Y,Z)
       =(R^\nabla(U,X)\cdot K)(Y,Z),
\end{equation}
where
\begin{equation}\label{eq:R-nabla_K}
    \begin{aligned}
        (R^\nabla(U,X)\cdot K)(Y,Z)
         ={}&R^\nabla(U,X)K(Y,Z)-K(R^\nabla(U,X)Y,Z)\\
            &-K(Y,R^\nabla(U,X)Z).
    \end{aligned}    
\end{equation}

\section{Hessian sectional curvature on conjugate symmetric statistical manifolds}
\label{sec:curvature}

In this section, we introduce Hessian sectional curvature,
an invariant originally defined by Shima~\cite{Shima_1995} for Hessian
manifolds, in the more general setting of conjugate symmetric
statistical manifolds.
We also define principal Hessian curvatures.

\subsection{Hessian curvature tensor}\label{subsec:Hessian_curvature_tensor}

To define Hessian sectional curvature, we first define the Hessian
curvature tensor and determine its symmetries.
For a statistical manifold $(M,g,\nabla)$, set
\begin{equation}\label{eq:Q_definition}
 Q(X,Y)Z:=-(\nabla_XK)(Y,Z).
\end{equation}
For Hessian manifolds, $Q$ is the \emph{Hessian curvature tensor} introduced by Shima~\cite{Shima}. 
We use the same terminology even when $\nabla$ is not flat.
We also use $Q$ to denote the associated $(0,4)$-tensor field
defined by
\begin{equation}\label{eq:Q_covariant}
 Q(W,X,Y,Z):=g(W, Q(X,Y)Z).
\end{equation}
Recall from~\eqref{eq:K-C} that
$\Kflat(X,Y,Z)=g(K(X,Y),Z)=-C(X,Y,Z)/2$.

\begin{Prop}\label{prop:Q_expansion}
On a statistical manifold $(M,g,\nabla)$,
\begin{align}
 Q(W,X,Y,Z)
 &=-(\nabla^g_X\Kflat)(W,Y,Z)+B_K(W,X,Y,Z),
 \label{eq:Q_expansion}\\
 \begin{split}
 B_K(W,X,Y,Z)
 &:=g(K(W,Y),K(X,Z))+g(K(W,Z),K(X,Y))\\
 &\quad-g(K(W,X),K(Y,Z)).
 \end{split}
 \label{eq:B_K}
\end{align}
Moreover, the following conditions are equivalent:
\begin{enumerate}[label=\textup{(\roman*)}]
 \item\label{item:Q_CS} $(M,g,\nabla)$ is conjugate symmetric;
 \item\label{item:Q_first_symmetry}
 $Q(W,X,Y,Z)=Q(X,W,Y,Z)$ for all vector fields;
 \item\label{item:Q_all_symmetries} The $(0,4)$-tensor field $Q$ has the symmetries
 \begin{equation}\label{eq:Q_symmetries}
 Q(W,X,Y,Z)=Q(X,W,Y,Z)
 =Q(W,X,Z,Y)=Q(Y,Z,W,X).
 \end{equation}
\end{enumerate}
\end{Prop}

For the proof, we prepare the following equality, obtained from $\nabla = \nabla^g + K$:
\begin{equation}\label{eq:nabla_K_expansion}
\begin{split}
 (\nabla_XK)(Y,Z)
 ={}&(\nabla^g_XK)(Y,Z)+K(X,K(Y,Z))\\
    &-K(K(X,Y),Z)-K(Y,K(X,Z)).
\end{split}
\end{equation}

\begin{proof}[Proof of Proposition~\ref{prop:Q_expansion}]
Taking the $g$-inner product of~\eqref{eq:nabla_K_expansion} with $-W$
and using the symmetry of $\Kflat$, we obtain~\eqref{eq:Q_expansion}.
$B_K$ has all the symmetries in~\eqref{eq:Q_symmetries}.
Since $\nabla^g_X\Kflat$ is symmetric in its three tensor arguments,
condition~\ref{item:Q_first_symmetry} is equivalent to total symmetry of
$\nabla^g\Kflat$, hence to conjugate symmetry by~\eqref{eq:CS_Codazzi}.
Total symmetry of $\nabla^g\Kflat$ also implies
condition~\ref{item:Q_all_symmetries}.
Clearly, condition~\ref{item:Q_all_symmetries} implies
condition~\ref{item:Q_first_symmetry}. This completes the proof.
\end{proof}

For a conjugate symmetric structure, the definition $Q=-\nabla K$,
together with \eqref{eq:CS_Codazzi} and \eqref{eq:nabla_K_expansion}, yields 
\begin{equation}\label{eq:Q_K_commutator}
 Q(X,Y)Z-Q(Y,X)Z=-2[K_X,K_Y]Z.
\end{equation}

\subsection{Hessian sectional curvature and principal Hessian curvatures}
\label{subsec:Hessian_sectional_principal}

In this subsection, we introduce Hessian sectional curvature and
principal Hessian curvatures for conjugate symmetric statistical manifolds.
Throughout this subsection, we fix a conjugate symmetric statistical
manifold $(M^n,g,\nabla)$.

Recall that we denote by $\LR_g$ the fiber metric on $S^2TM$
induced by $g$. In particular,
\begin{equation}\label{eq:symmetric_metric}
 \langle W\odot X,Y\odot Z\rangle_g
 =\frac12\{g(W,Y)g(X,Z)+g(W,Z)g(X,Y)\}.
\end{equation}
By Proposition~\ref{prop:Q_expansion}, the formula
\[
 \calQ(W\odot X,Y\odot Z):=Q(W,X,Y,Z)
\]
defines a symmetric bilinear form on $S^2TM$.
In this paper, we call the self-adjoint bundle endomorphism ${\hat{Q}}^{(g,\nabla)}:S^2TM\to S^2TM$ characterized by
\begin{equation}\label{eq:A_definition}
 \langle {\hat{Q}}^{(g,\nabla)} (\nu), \mu \rangle_g = \calQ(\nu,\mu)
\end{equation}
for all $\nu,\mu\in\Gamma(S^2TM)$
the \emph{Hessian curvature operator}.
We also write $\hat{Q}$ when the structure is clear.

\begin{Def}\label{def:Hessian_curvatures}
For $p\in M$ and $0\ne\nu\in S^2T_pM$, the \emph{Hessian sectional
curvature in the direction $\nu$} is
\begin{equation}\label{eq:Hessian_quotient}
 q_p(\nu):=\frac{\calQ_p(\nu,\nu)}{\langle \nu,\nu \rangle_{g_p}} = \frac{\langle \hat{Q}_p (\nu), \nu \rangle_{g_p}}{\langle \nu,\nu \rangle_{g_p}} .
\end{equation}
We say that each eigenvalue of $\hat{Q}_p$ is a
\emph{principal Hessian curvature at $p$}.
The multiset of these eigenvalues, counted with multiplicity, is denoted
by $\PHC_p$, or by $\PHC_p(g,\nabla)$ when the structure is specified.
\end{Def}

For each $p\in M$, the value $q_p(\nu)$ depends only on the line
spanned by $\nu$, and hence $q$ defines a smooth function on the
projective bundle $\mathbb P(S^2TM)$.

In the Hessian case, the Hessian sectional curvature defined above
coincides with that introduced by Shima~\cite{Shima_1995}.
Principal Hessian curvatures are also essentially the same as
the \emph{principal conical curvatures} introduced by Vaisman
in~\cite[Open Problem in Appendix]{Vaisman_2014}.
However, owing to different conventions for the fiber metric on
$S^2TM$, our principal Hessian curvatures are twice the corresponding principal conical curvatures.

The following holds. The proof is straightforward, so we omit it.

\begin{Thm}\label{thm:HSC-def-on-CS}
Both the Hessian sectional curvature and the principal Hessian
curvatures are invariant under statistical isomorphisms.
\end{Thm}

We now describe the relation between Hessian sectional curvature
and principal Hessian curvatures.
Let $\{v_1,\ldots,v_N\}$ be an orthonormal eigenbasis of
$\hat{Q}_p$, where $N=n(n+1)/2$, with corresponding
eigenvalues $\kappa_1,\ldots,\kappa_N$.
Then
\begin{equation}\label{eq:Rayleigh}
 q_p\left(\sum_{j=1}^N t_jv_j\right)
 =\frac{\sum_{j=1}^N\kappa_jt_j^2}{\sum_{j=1}^N t_j^2}
 \qquad\bigl((t_1,\ldots,t_N)\ne0\bigr).
\end{equation}
Thus, the principal Hessian curvatures together with their
corresponding eigenspaces determine $q_p$.
Conversely, $q_p$ determines $\calQ_p$ by polarization and hence
determines $\hat{Q}_p$.

\subsection{Multiplicities of principal Hessian curvatures}
\label{subsec:parallel_cubic_curvature}

In this subsection, we study the multiplicities of principal
Hessian curvatures.

Let $(M,g,\nabla)$ be a statistical manifold with difference tensor $K$.
Define the bundle map
\begin{equation}\label{eq:tilcalK_def}
 \calK:S^2TM\longrightarrow TM,\qquad
 \calK(X\odot Y):=K(X,Y).
\end{equation}
We denote its adjoint with respect to $\LR_g$ and $g$ by
$\calK^*:TM\to S^2TM$.
For any local $g$-orthonormal frame $\{e_1,\ldots,e_n\}$, set
\begin{equation}\label{eq:trace_symmetric_tensor}
 g^\sharp:=\sum_{i=1}^n e_i\odot e_i,\qquad
 \tr_g\nu:=\langle \nu,g^\sharp\rangle_g
 \quad (\nu\in\Gamma(S^2TM)).
\end{equation}

\begin{Lem}\label{lem:gram_curvature}
Let $(M,g,\nabla)$ be a statistical manifold with $\nabla^g C=0$,
which is therefore conjugate symmetric (see Section~\ref{subsec:class_structures}).
Suppose it has pointwise constant sectional $\mathsf K$-curvature
$\rho$, that is,
\begin{equation}\label{eq:gram_commutator_assumption}
 [K_X,K_Y]Z=\rho\{g(Y,Z)X-g(X,Z)Y\}.
\end{equation}
Then
\begin{equation}\label{eq:gram_curvature}
 \calQ(\nu,\mu)
 =2\rho\{ \langle \nu,\mu \rangle_g -(\tr_g\nu)(\tr_g\mu)\}
   +g(\calK\nu,\calK\mu).
\end{equation}
Equivalently,
\begin{equation}\label{eq:gram_operator}
 \hat{Q} (\nu) =2\rho\{\nu-(\tr_g\nu)g^\sharp\}
       +(\calK^* \circ \calK)(\nu).
\end{equation}
\end{Lem}

\begin{proof}
Taking the $g$-inner product of
\eqref{eq:gram_commutator_assumption} with $W$ gives
\begin{align*}
 &g(K(W,X),K(Y,Z))-g(K(W,Y),K(X,Z))\\
 &\hspace{10mm}=\rho\{g(W,X)g(Y,Z)-g(W,Y)g(X,Z)\}.
\end{align*}
Interchanging $Y$ and $Z$ in the preceding identity and using
the symmetry of $K$ and $g$, we obtain
\begin{align*}
 &g(K(W,X),K(Y,Z))-g(K(W,Z),K(X,Y))\\
 &\hspace{10mm}
 =\rho\{g(W,X)g(Y,Z)-g(W,Z)g(X,Y)\}.
\end{align*}
Substituting these two identities into~\eqref{eq:B_K} and using
$\nabla^g\Kflat=0$, we obtain from~\eqref{eq:Q_expansion}
\begin{align*}
 Q(W,X,Y,Z)
 ={}&\rho\{g(W,Y)g(X,Z)+g(W,Z)g(X,Y)\\
 &\hspace{19mm}-2g(W,X)g(Y,Z)\}\\
 &+g(K(W,X),K(Y,Z)).
\end{align*}
Now take $\nu=W\odot X$, $\mu=Y\odot Z$ and use
\eqref{eq:symmetric_metric},
$\tr_g(W\odot X)=g(W,X)$, and~\eqref{eq:tilcalK_def}.
The resulting identity is~\eqref{eq:gram_curvature}.
By bilinearity, this yields~\eqref{eq:gram_curvature} for arbitrary $\nu,\mu\in \Gamma(S^2TM)$.
\end{proof}

The following is an immediate consequence of~\eqref{eq:gram_operator} and the rank--nullity theorem.

\begin{Cor}\label{cor:PHC_rank_bound}
Under the assumptions of Lemma~\ref{lem:gram_curvature},
if $n\geq3$, then $2\rho(p)$ is a principal Hessian curvature
at $p$ with multiplicity at least $(n-2)(n+1)/2$.
If $\rho(p)=0$ and $n\geq2$, then $0$ is a principal Hessian
curvature at $p$ with multiplicity at least $n(n-1)/2$.
\end{Cor}

The explicit computation in Theorem~\ref{thm:euclidean_spectrum} shows
that, for $n\geq2$, every left-invariant Hessian structure on $\RR^n$
has the principal Hessian curvature $0$ with multiplicity at least
$n(n-1)/2$. Corollary~\ref{cor:PHC_rank_bound} explains this as a consequence of
$\nabla^g C=0$ and vanishing sectional $\mathsf K$-curvature, that is $[K, K] = 0$.

We conclude this section with an example showing that principal Hessian
curvatures can distinguish non-Hessian statistical structures that share
several other invariants.

\begin{Ex}\label{ex:PHC_distinguishes}
Let $g_E=(dx^1)^2+(dx^2)^2$ be the Euclidean metric on $\RR^2$, and
let $e_i=\partial/\partial x^i$ be the standard frame.
For constants $u,v\in\RR$, define a symmetric $(1,2)$-tensor field $K$ by
\[
 K(e_1,e_1)=u e_1,\qquad K(e_1,e_2)=v e_2,\qquad
 K(e_2,e_2)=v e_1.
\]
Then $\nabla:=\nabla^{g_E}+K$ defines a conjugate symmetric
statistical structure satisfying
\[
 \nabla^{g_E}K=0,\qquad R^\nabla=[K,K]=v(u-v)R_0.
\]
Consider $(u,v)=(0,2)$ and $(u,v)=(3,-1)$.
In both cases,
\[
 R^\nabla=-4R_0,\qquad \lVert C \rVert_{g_E}^2=48,\qquad
 \tr_{g_E}K=2e_1.
\]
By Lemma~\ref{lem:gram_curvature}, the matrices of their Hessian
curvature operators in the orthonormal frame $e_1\odot e_1,\sqrt2e_1\odot e_2,e_2\odot e_2$ are
\[
 A_1=\begin{pmatrix}0&0&8\\0&0&0\\8&0&4\end{pmatrix},
 \qquad
 A_2=\begin{pmatrix}9&0&5\\0&-6&0\\5&0&1\end{pmatrix}.
\]
The first has a zero eigenvalue, whereas $\det A_2=96\ne0$.
The two structures are therefore not statistically isomorphic.
\end{Ex}

\section{Constant and weakly constant Hessian sectional curvature}
\label{sec:conditions}

In this section, we introduce and study the notions of constant
Hessian sectional curvature (CHC) and weakly constant Hessian
sectional curvature (weakly CHC) in the setting of conjugate
symmetric statistical manifolds.

\subsection{Definitions and characterizations}\label{subsec:CHC_characterizations}

In this subsection, we give the definitions and characterizations
of these notions.

\begin{Def}\label{def:CHC}
A conjugate symmetric statistical manifold $(M,g,\nabla)$ is 
\emph{of constant Hessian sectional curvature $c$}
(\emph{of CHC $c$} for short) if $q_p(\nu)$ takes the same value
$c\in\RR$ for every $p\in M$ and every nonzero $\nu\in S^2T_pM$.
For $p\in M$, we say that $(M,g,\nabla)$ is \emph{of CHC at $p$}
if $q_p$ is constant.
It is \emph{pointwise CHC} if it is of CHC at every point.
\end{Def}

For Hessian manifolds, this agrees with the definition
by Shima~\cite{Shima_1995}.
Examples of Hessian manifolds of CHC can be found, for instance,
in~\cite{Furuhata_2011,FK_2013,Shima_1995}.
With our extension to conjugate symmetric statistical manifolds, every
trivial statistical structure is of CHC $0$, since $K=0$.
Non-trivial statistical structures of CHC, including non-Hessian examples, will be constructed in Section~\ref{subsec:non-Hessian_Ex}.
A Schur-type theorem for Hessian sectional curvature is given in Proposition~\ref{prop:CHC_rigidity}~\ref{item:ptCHC-CHC_rigidity}.

The following characterization follows from~\eqref{eq:symmetric_metric},
\eqref{eq:A_definition}, and~\eqref{eq:Rayleigh}.

\begin{Prop}\label{prop:CHC_tensor}
Let $(M,g,\nabla)$ be a conjugate symmetric statistical manifold and let
$c\in\RR$. The following conditions are equivalent:
\begin{enumerate}[label=\textup{(\roman*)}]
 \item\label{item:CHC_quotient} $(M,g,\nabla)$ is of CHC $c$;
 \item\label{item:CHC_operator} $\hat{Q} = c\,\id$;
 \item\label{item:CHC_tensor} for all vector fields $X,Y,Z$,
 \begin{equation}\label{eq:CHC_tensor}
 Q(X,Y)Z=\frac c2\{g(X,Y)Z+g(X,Z)Y\}.
 \end{equation}
\end{enumerate}
The same equivalences hold pointwise, with $c$ replaced by a smooth
function on $M$.
\end{Prop}

In the Hessian case, Shima~\cite[Theorem~1]{Shima_1995} proved
the equivalence of~\ref{item:CHC_quotient} and~\ref{item:CHC_tensor}.
The pointwise form of~\ref{item:CHC_operator} appears explicitly
in the proof.

By Proposition~\ref{prop:CHC_tensor}, the following holds.

\begin{Cor}\label{cor:CHC_equal_PHC}
A conjugate symmetric statistical manifold is of CHC at a point $p$ if
and only if all principal Hessian curvatures at $p$ are equal.
In particular, it is pointwise CHC if and only if this equality holds
at every point.
\end{Cor}

If a conjugate symmetric statistical manifold $(M^n,g,\nabla)$ with
$n\geq2$ is pointwise CHC, then the function $c$ in
Proposition~\ref{prop:CHC_tensor} is given by
\begin{equation}\label{eq:CHC-c_Scals}
    c=\frac{4\bigl(\operatorname{Scal}^{\nabla}
               -\operatorname{Scal}^{g}\bigr)}{n(n-1)}.
\end{equation}
This follows by combining~\eqref{eq:conjugate_curvature},
\eqref{eq:Q_K_commutator}, and~\eqref{eq:CHC_tensor}
and taking traces.

We next introduce the weakly CHC condition.
For an orthogonal map $u:T_pM\to T_pM$, denote its induced action on $S^2T_pM$ by $u\cdot(X\odot Y)=uX\odot uY$.

\begin{Def}\label{def:weak_CHC}
A conjugate symmetric statistical manifold $(M,g,\nabla)$ has
\emph{weakly constant Hessian sectional curvature at $p\in M$}
(or is \emph{weakly CHC at $p$}) if
\[
 q_p(u\cdot\nu)=q_p(\nu)
 \qquad\bigl(u\in O(T_pM,g_p),\ 0\ne\nu\in S^2T_pM\bigr).
\]
It is \emph{weakly CHC} if this holds at every point.
\end{Def}

Every pointwise CHC conjugate symmetric statistical manifold is weakly CHC.

\begin{Rem}\label{rem:dimension_one}
Every $1$-dimensional statistical manifold is Hessian.
It is also pointwise CHC, since $S^2T_pM$ is $1$-dimensional
at every point $p$, and hence weakly CHC.
However, the Hessian sectional curvature need not be constant even
when $M$ is connected. 
For example, equip $\RR$ with its standard flat connection $D$
and the metric $g=e^{x^2}(dx)^2$.
By~\cite[Proposition~1]{Molitor_2022}, the Hessian sectional
curvature at $x$ is $e^{-x^2}$.
Thus the Hessian manifold $(\RR, g, D)$ is not of CHC.
\end{Rem}

For $n\geq2$, write
\[
 S^2TM=\RR g^\sharp\oplus S^2_0TM,\qquad S^2_0TM:=\ker(\tr_g),
\]
where $g^\sharp$ and $\tr_g$ are defined
in~\eqref{eq:trace_symmetric_tensor}, and let
\[
 \pr_{\tr}(\nu):=\frac{\tr_g\nu}{n}g^\sharp,
 \qquad \pr_0:=\id-\pr_{\tr}
\]
be the corresponding orthogonal projections.

\begin{Prop}\label{prop:weak_characterization}
Let $(M^n,g,\nabla)$ be conjugate symmetric, with $n\geq2$, and let
$\hat{Q}$ be its Hessian curvature operator, defined by~\eqref{eq:A_definition}.
The following conditions are equivalent:
\begin{enumerate}[label=\textup{(\roman*)}]
 \item\label{item:weak_quotient} $(M,g,\nabla)$ is weakly CHC;
 \item\label{item:weak_equivariant}
 $\hat{Q}_p$ is $O(T_pM,g_p)$-equivariant for every $p$;
 \item\label{item:weak_operator}
 there are unique smooth functions $\lambda_{\tr},\lambda_0$ on $M$ such that
 \begin{equation}\label{eq:weak_operator}
 \hat{Q} = \lambda_{\tr}\pr_{\tr}+\lambda_0\pr_0;
 \end{equation}
 \item\label{item:weak_tensor}
 there are unique smooth functions $a,b$ such that
 \begin{equation}\label{eq:weak_tensor}
 Q(X,Y)Z=a\{g(X,Y)Z+g(X,Z)Y\}+b\,g(Y,Z)X.
 \end{equation}
\end{enumerate}
The coefficients are related by
\begin{equation}\label{eq:weak_coefficients}
 \lambda_0=2a,\qquad \lambda_{\tr}=2a+nb.
\end{equation}
The structure is pointwise CHC exactly when $b=0$, and it is of CHC $c$
exactly when $b=0$ and $2a\equiv c$.
\end{Prop}

\begin{proof}
The equivalence
\ref{item:weak_quotient}$\Leftrightarrow$\ref{item:weak_equivariant}
follows from polarization of the quadratic form $\nu\mapsto\calQ_p(\nu,\nu)$ and the orthogonal invariance
of $\LR_{g_p}$.

The decomposition
$S^2T_pM=\RR g_p^\sharp\oplus S^2_0T_pM$
is a direct sum of two inequivalent irreducible real representations
of $O(T_pM,g_p)$.
Since $\hat{Q}_p$ is self-adjoint, Schur's lemma (see, for example,~\cite[Chapter~III, Section~1]{Weyl_1939}) implies that \ref{item:weak_equivariant} is equivalent, pointwise,
to \ref{item:weak_operator}, or equivalently to
\begin{equation}\label{eq:weak_bilinear}
 \calQ(\nu,\mu)
 =\lambda_0 \langle \nu,\mu \rangle_g
  +\frac{\lambda_{\tr}-\lambda_0}{n}(\tr_g\nu)(\tr_g\mu).
\end{equation}
The coefficients are unique because both summands have positive
dimension for $n\geq2$, and they are smooth by the formulas
\begin{equation}
    \lambda_{\tr}=\frac1n \langle \hat{Q} (g^\sharp) , g^\sharp \rangle_g ,\qquad \lambda_0=\frac{\tr \hat{Q} -\lambda_{\tr}}{n(n+1)/2-1}.
\end{equation}
Evaluating~\eqref{eq:weak_bilinear} on $W\odot X$ and $Y\odot Z$
and using~\eqref{eq:symmetric_metric} yields
\ref{item:weak_operator}$\Leftrightarrow$\ref{item:weak_tensor}
and~\eqref{eq:weak_coefficients}.
The last assertions follow from Proposition~\ref{prop:CHC_tensor}.
\end{proof}

The tensor identity~\eqref{eq:weak_tensor} also makes sense without
assuming conjugate symmetry. 
By Proposition~\ref{prop:Q_expansion}, specifically the equivalence of
\ref{item:Q_CS} and \ref{item:Q_first_symmetry}, this identity itself implies conjugate symmetry.

\begin{Cor}\label{cor:weak_forces_CS}
Let $(M^n,g,\nabla)$ be any statistical manifold, with $n\geq2$.
If its Hessian curvature tensor $Q$ satisfies~\eqref{eq:weak_tensor}
for some functions $a,b$, then it is conjugate symmetric and
weakly CHC.
\end{Cor}

We conclude this subsection with the following lemma, which will be used in the proofs in the next subsection.

\begin{Lem}\label{lem:weak_integrability}
Let $(M^n,g,\nabla)$ be a weakly CHC conjugate symmetric statistical
manifold with $n\geq2$.
Let $a,b$ be the functions in Proposition~\ref{prop:weak_characterization}.
Then
\begin{equation}\label{eq:weak_integrability}
\begin{aligned}
 -\bigl(R^\nabla(U,X)\cdot K\bigr)(Y,Z)
 ={}&da(U)\{g(X,Y)Z+g(X,Z)Y\}\\
 &-da(X)\{g(U,Y)Z+g(U,Z)Y\}\\
 &+\{db(U)g(Y,Z)-2b\Kflat(U,Y,Z)\}X\\
 &-\{db(X)g(Y,Z)-2b\Kflat(X,Y,Z)\}U
\end{aligned}
\end{equation}
for all vector fields $U,X,Y,Z$.
In particular, if $(M,g,\nabla)$ is of CHC, then
\begin{equation}\label{eq:CHC_curvature_K}
 R^\nabla(U,X)\cdot K=0
\end{equation}
for all vector fields $U,X$ (see~\eqref{eq:R-nabla_K} for the notation $R^\nabla(U,X)\cdot K$).
\end{Lem}

\begin{proof}
Since $Q=-\nabla K$, the Ricci identity
\eqref{eq:Ricci_commutation} gives
\[
 (\nabla_UQ)(X,Y)Z-(\nabla_XQ)(U,Y)Z
 =-\bigl(R^\nabla(U,X)\cdot K\bigr)(Y,Z).
\]
Substituting~\eqref{eq:weak_tensor} and using
$\nabla g=-2\Kflat$ and the total symmetry of $\Kflat$,
we obtain~\eqref{eq:weak_integrability}.
If the structure is of CHC $c$, then
$a=c/2$ and $b=0$ by
Proposition~\ref{prop:weak_characterization}.
Since $c$ is constant,
equation~\eqref{eq:weak_integrability}
reduces to~\eqref{eq:CHC_curvature_K}.
\end{proof}

\subsection{CHC conditions and other notions of constant curvature}\label{subsec:K_curvature_Schur}

In this subsection, we first study the relations of the (weakly) CHC
condition to constant sectional $\mathsf K$-curvature and constant
curvature.
We then investigate the relations among the weakly CHC, pointwise CHC,
and CHC conditions (see~Proposition~\ref{prop:CHC_rigidity}).

We begin with the following lemma relating the weakly CHC condition
to sectional $\mathsf K$-curvature.

\begin{Lem}\label{lem:weak_K_curvature}
Let $(M^n,g,\nabla)$ be a weakly CHC conjugate symmetric statistical manifold
with $n\geq2$. Then
\begin{equation}\label{eq:weak_rho}
 [K_X,K_Y]=\rho R_0(X,Y),\qquad
 \rho := \frac{a-b}{2}
 =\frac{(n+2)\lambda_0-2\lambda_{\tr}}{4n}.
\end{equation}
Here $a,b,\lambda_0,\lambda_{\tr}$ are the functions in
Proposition~\ref{prop:weak_characterization}.
If $M$ is connected and $n\geq3$, then $\rho$ is constant.
\end{Lem}

\begin{proof}
The first identity follows from~\eqref{eq:Q_K_commutator} and~\eqref{eq:weak_tensor}.
The equality of the two expressions for $\rho$ follows from~\eqref{eq:weak_coefficients}.
The last assertion follows from Proposition~\ref{prop:Opozda_K_Schur}.
\end{proof}

Suppose that $(M,g,\nabla)$ is conjugate symmetric and of CHC $c$.
Proposition~\ref{prop:weak_characterization} and Lemma~\ref{lem:weak_K_curvature} give 
$$
    [K,K] = \frac{c}{4} R_0.
$$
Together with~\eqref{eq:conjugate_curvature}, we have the following.

\begin{Prop}\label{prop:CHC_curvature_relation}
Let $(M^n,g,\nabla)$ be a conjugate symmetric statistical
manifold with $n\geq2$.
Suppose that the structure is of CHC $c$.
Then
\begin{equation}\label{eq:CHC_curvature_relation}
    [K,K] = R^\nabla-R^g=\frac c4R_0.
\end{equation}
In particular, $(M,g,\nabla)$ is Hessian if and only if
$g$ has constant sectional curvature $-c/4$.
\end{Prop}

Note that it is known from Shima~\cite[Corollary~2]{Shima_1995} that, in the Hessian case, CHC $c$ implies that $g$ has constant sectional curvature $-c/4$.

We next consider the case where $R^\nabla=kR_0$
for a nonzero constant $k$.

\begin{Prop}\label{prop:CHC_nonzero_affine_curvature}
Let $(M^n,g,\nabla)$ be conjugate symmetric, with $n\geq2$.
If it is of CHC $c$ and has constant curvature $k\ne0$, that is,
$R^\nabla=kR_0$, then $K=0$ and $c=0$.
\end{Prop}

\begin{proof}
By Lemma~\ref{lem:weak_integrability}, $R^\nabla(U,X)\cdot K=0$.
Since $R^\nabla(U,X)$ is $g$-skew-adjoint (see~Section~\ref{subsec:class_structures}), this implies $R^\nabla(U,X)\cdot\Kflat=0$.
Fix $p\in M$.
As $k\ne0$, the endomorphisms $k (R_0)_p (U,X)$ span
$\mathfrak{so}(T_pM,g_p)$, thus $B\cdot\Kflat_p=0$ for every
$B\in\mathfrak{so}(T_pM,g_p)$.
Since $SO(T_pM,g_p)$ is connected,
$\Kflat_p$ is $SO(T_pM,g_p)$-invariant
(cf.~\cite[Section~4.2]{Kirillov_2008}).
Since $n\geq2$, this group acts transitively on the unit sphere in $T_pM$.
It follows that the function
\[
 \Phi(x):=\Kflat_p(x,x,x)=g_p(K_p(x,x),x)
\]
is constant on that sphere.
Since $\Phi(-x)=-\Phi(x)$, this constant must be zero.
Thus $\Kflat_p=0$.
Since $p$ was arbitrary, $K=0$ on $M$ and hence $c=0$.
\end{proof}

The preceding two propositions give the following corollary.

\begin{Cor}\label{cor:CHC_nonHessian_sectional_curvature}
Let $(M^n,g,\nabla)$ be a non-trivial conjugate symmetric statistical
manifold of CHC $c$, with $n\geq2$.
If $(M,g,\nabla)$ is not Hessian, then $g$ does not have constant
sectional curvature.
\end{Cor}

\begin{proof}
Suppose that $g$ has constant sectional curvature $\kappa$.
Then~\eqref{eq:CHC_curvature_relation} gives
$R^\nabla=(\kappa+c/4)R_0$.
Since $(M,g,\nabla)$ is not Hessian, $\kappa+c/4\ne0$.
Proposition~\ref{prop:CHC_nonzero_affine_curvature} now gives $K=0$,
contrary to the assumption that the structure is non-trivial.
\end{proof}

The following proposition also concerns statistical structures of CHC.
It characterizes when the constant Hessian sectional curvature is negative. 

\begin{Prop}\label{prop:CHC_sign}
Let $(M^n,g,\nabla)$ be a conjugate symmetric statistical
manifold of CHC $c$, with $n\geq2$.
Then
\begin{equation}\label{eq:CHC_scalar_norm}
    \lVert\tr_gK\rVert_g^2-\lVert K\rVert_g^2 
    = \operatorname{Scal}^{\nabla}-\operatorname{Scal}^{g} 
    = \frac{c}{4}n(n -1).
\end{equation}
Moreover, the following conditions are equivalent:
\begin{enumerate}[label=\textup{(\roman*)}]
 \item\label{item:CHC_sign_negative}
 $c<0$;
 \item\label{item:CHC_sign_norm}
$|K_p|^2 > |(\tr_gK)_p|^2$
at some point $p\in M$;
 \item\label{item:CHC_sign_hessian_positive_scalar}
 $(M,g,\nabla)$ is Hessian and $\operatorname{Scal}^{g}>0$.
\end{enumerate}
In particular, if $(\tr_gK)_p=0$ and $K_p\ne0$
at some point $p\in M$, then $c<0$.
\end{Prop}

\begin{proof}
Equation~\eqref{eq:CHC_scalar_norm} follows from \eqref{eq:stat_scalar_norm} and~\eqref{eq:CHC-c_Scals}.
Since $c$ is constant and $n\geq2$, this identity gives the equivalence of
\ref{item:CHC_sign_negative} and \ref{item:CHC_sign_norm}.
The last assertion is a special case of the implication
\ref{item:CHC_sign_norm}~$\Rightarrow$~\ref{item:CHC_sign_negative}.

In the Hessian case, \eqref{eq:CHC_scalar_norm} gives
$\operatorname{Scal}^{g}=-n(n-1)c/4$, so \ref{item:CHC_sign_hessian_positive_scalar} implies \ref{item:CHC_sign_negative}.
Conversely, suppose that $c<0$.
By~\eqref{eq:CHC_curvature_relation}, the structure has constant
sectional $\mathsf K$-curvature $c/4<0$.
Each $R^\nabla(U,X)$ is $g$-skew-adjoint by conjugate symmetry
(see Section~\ref{subsec:class_structures}) and satisfies
$R^\nabla(U,X)\cdot K=0$ by~\eqref{eq:CHC_curvature_K}.
Proposition~\ref{prop:Opozda_negative_K} therefore gives $R^\nabla=0$.
Thus the structure is Hessian, and~\eqref{eq:CHC_scalar_norm}
gives $\operatorname{Scal}^{g}>0$.
\end{proof}

Combining Proposition~\ref{prop:CHC_sign}
with a result of Shima~\cite[Corollary~3]{Shima_1995},
we obtain the following corollary.

\begin{Cor}\label{cor:complete_CHC_nonnegative}
Let $(M^n,g,\nabla)$ be a conjugate symmetric statistical
manifold of CHC $c$, with $n\geq2$.
If $g$ is complete, then $c\geq0$.
\end{Cor}

\begin{proof}
Suppose that $c<0$.
By Proposition~\ref{prop:CHC_sign}, the structure is Hessian.
The universal cover of any connected component of $M$,
equipped with the lifted statistical structure, is a simply connected
Hessian manifold of CHC $c$.
Since $g$ is complete, so is the lifted Riemannian metric.
This contradicts~\cite[Corollary~3]{Shima_1995}.
\end{proof}

We conclude this subsection with Proposition~\ref{prop:CHC_rigidity},
which describes the relations among the weakly CHC, pointwise CHC, and CHC conditions. For the proof, we prepare the following lemma.

\begin{Lem}\label{lem:weak_CHC_identities}
Let $(M^n,g,\nabla)$ be a weakly CHC conjugate symmetric
statistical manifold with $n\geq2$.
Let $a,b$ be the functions in
Proposition~\ref{prop:weak_characterization},
and put $\rho:=(a-b)/2$.
Define the $1$-form $\tau$ by $\tau(U):=\tr K_U$.
Then
\begin{align}
 (n-1)d\rho&=-b\tau,
 \label{eq:weak_CHC_drho}\\
 (n-2)bK&=0.
 \label{eq:weak_CHC_bK}
\end{align}
Moreover,
\begin{equation}\label{eq:weak_CHC_nablaK}
 (\nabla^g_XK)(Y,Z)
 =K(K(X,Y),Z)-a\,g(X,Y)Z - \frac{a+b}{2}\{g(X,Z)Y+g(Y,Z)X\}.
\end{equation}
\end{Lem}

\begin{proof}
Equation~\eqref{eq:weak_CHC_nablaK} follows from \eqref{eq:nabla_K_expansion}, \eqref{eq:weak_tensor}, and~\eqref{eq:weak_rho}.
Differentiating $[K_X,K_Y]=\rho R_0(X,Y)$ with respect to
$\nabla^g_U$, and using $\nabla^gR_0=0$, \eqref{eq:weak_CHC_nablaK}, and the total symmetry of $\Kflat$, we obtain
\begin{equation}\label{eq:weak_rigidity_commutator}
\begin{aligned}
 d\rho(U)R_0(X,Y)
 &=[(\nabla^g_UK)_X,K_Y]+[K_X,(\nabla^g_UK)_Y]\\
 &=-b\{R_0(K(U,X),Y)+R_0(X,K(U,Y))\}.
\end{aligned}
\end{equation}
Applying this identity to $Z$ and taking the trace of the
resulting endomorphism in $X$ yields
\begin{equation}\label{eq:weak_rigidity_trace}
 \bigl((n-1)d\rho(U)+b\tau(U)\bigr)g(Y,Z)
 =-(n-2)b\Kflat(U,Y,Z).
\end{equation}
Put $\eta:=(n-1)d\rho+b\tau$.
Since $\Kflat$ is totally symmetric, interchanging $U$ and $Y$
in~\eqref{eq:weak_rigidity_trace} gives
\[
 \eta(U)g(Y,Z)=\eta(Y)g(U,Z).
\]
Taking the $g$-trace in $Y$ and $Z$, we obtain
$(n-1)\eta(U)=0$.
Since $n\geq2$, we have $\eta=0$, proving
\eqref{eq:weak_CHC_drho}.
Substitution into~\eqref{eq:weak_rigidity_trace}
gives~\eqref{eq:weak_CHC_bK}.
\end{proof}

\begin{Prop}\label{prop:CHC_rigidity}
Let $(M^n,g,\nabla)$ be a connected conjugate symmetric
statistical manifold with $n\geq2$.
Then the following statements hold.
\begin{enumerate}[label=\textup{(\arabic*)}]
 \item\label{item:CHC_rigidity_wCHC-CHC}
 If $n\geq3$, then $(M,g,\nabla)$ is weakly CHC if and only if it is of CHC.
 \item \label{item:ptCHC-CHC_rigidity}
 $(M,g,\nabla)$ is pointwise CHC if and only if it is of CHC.
\end{enumerate}
\end{Prop}

For connected Hessian manifolds of dimension at least $2$, Shima
proved that pointwise CHC implies CHC \cite[Theorem~1]{Shima_1995}.
Proposition~\ref{prop:CHC_rigidity}~\ref{item:ptCHC-CHC_rigidity} extends this result to connected conjugate symmetric statistical manifolds.

We now give the proof.

\begin{proof}[Proof of Proposition~\ref{prop:CHC_rigidity}]
It suffices to prove the ``only if'' part of each statement.
Suppose that $(M,g,\nabla)$ is weakly CHC, and use the notation
of Lemma~\ref{lem:weak_CHC_identities}.

We first prove assertion~\ref{item:ptCHC-CHC_rigidity}.
If the structure is pointwise CHC, then $b=0$ by
Proposition~\ref{prop:weak_characterization}.
Equation~\eqref{eq:weak_CHC_drho} gives $d\rho=0$.
Since $M$ is connected, $c:=4\rho=2a$ is constant.
Thus Proposition~\ref{prop:weak_characterization}
shows that the structure is of CHC $c$.

To prove assertion~\ref{item:CHC_rigidity_wCHC-CHC},
suppose that $n\geq3$.
Equation~\eqref{eq:weak_CHC_bK} gives $bK=0$.
On the open set $\Omega:=\{p\in M\mid b(p)\ne0\}$,
we therefore have $K=0$, and hence $Q=-\nabla K=0$.
The uniqueness of the coefficients in
Proposition~\ref{prop:weak_characterization}
implies that $a=b=0$ on $\Omega$.
Since $b$ is nonzero on $\Omega$ by definition,
we must have $\Omega=\emptyset$.
Thus $b=0$ everywhere, so the structure is pointwise CHC
by Proposition~\ref{prop:weak_characterization}.
Assertion~\ref{item:ptCHC-CHC_rigidity} now shows that it is of CHC.
\end{proof}

For connected Hessian manifolds of dimension at least $2$, the weakly CHC
and CHC conditions are equivalent (see~Theorem~\ref{thm:hessian_weak}).
In dimension $2$, however, there exist connected non-Hessian
conjugate symmetric statistical manifolds that are weakly CHC
but are not of CHC. 
Such examples are given in Section~\ref{sec:parallel}.

\subsection{The Hessian case}\label{subsec:Hessian_weak_CHC}

In this subsection, we show that CHC and weakly CHC coincide for connected Hessian
manifolds of dimension at least $2$ (Theorem~\ref{thm:hessian_weak}).
For its proof, we prepare the following lemma.

\begin{Lem}\label{lem:pure_trace_vanishing}
Let $(M^n,g,\nabla)$ be a weakly CHC conjugate symmetric
statistical manifold with $n\geq2$, and let $U\subset M$ be open.
If on $U$ the tensor $\Kflat$ has the form
\begin{equation}\label{eq:pure_trace}
 \Kflat(X,Y,Z)
 =\theta(X)g(Y,Z)+\theta(Y)g(X,Z)+\theta(Z)g(X,Y)
\end{equation}
for a $1$-form $\theta$ on $U$, then $K=0$ on $U$.
\end{Lem}

\begin{proof}
By Proposition~\ref{prop:pure_trace_Codazzi}, there exists a function $\sigma$ on $U$ such that $\nabla^g\theta=\sigma g$.
Let $E=\theta^\sharp$ be the vector field on $U$ defined by $g(E,X)=\theta(X)$ for all vector fields $X$ on $U$.
By~\eqref{eq:pure_trace}, we have
\begin{equation}\label{eq:pure_trace_K}
 K(X,Y)=\theta(X)Y+\theta(Y)X+g(X,Y)E.
\end{equation}
Now fix $p\in U$ and let $X\in T_pM$ be a unit vector.
By~\eqref{eq:pure_trace} and $\nabla^g\theta=\sigma g$, we have
\[
 (\nabla^g_X\Kflat)(X,X,X)=3\sigma(p).
\]
Moreover, by~\eqref{eq:pure_trace_K} and $g_p(E_p,X)=\theta_p(X)$,
\[
 |K_p(X,X)|^2
 =|2\theta_p(X)X+E_p|^2
 =|E_p|^2+8\theta_p(X)^2.
\]
By~\eqref{eq:Q_expansion}, we obtain
\[
 Q_p(X,X,X,X)
 =-(\nabla^g_X\Kflat)_p (X,X,X)+|K_p(X,X)|^2
 =-3\sigma(p)+|E_p|^2+8\theta_p(X)^2.
\]
Let $a,b$ be the functions in
Proposition~\ref{prop:weak_characterization}.
Then $Q_p(X,X,X,X)=2a(p)+b(p)$ for every unit vector $X\in T_pM$.
We show that $E_p=0$. 
Suppose, for a contradiction, that $E_p\ne0$.
Since $n\geq2$, we may choose unit vectors
$X=E_p/|E_p|$ and $Y\in T_pM$ with $Y\perp E_p$.
For these vectors, $\theta_p(X)^2=|E_p|^2$ and $\theta_p(Y)^2=0$.
Therefore,
\[
 Q_p(X,X,X,X)-Q_p(Y,Y,Y,Y)=8|E_p|^2>0.
\]
This contradicts the fact that both terms on the left-hand side
are equal to $2a(p)+b(p)$.
Thus $E_p=0$.
Since $p\in U$ was arbitrary, $E=0$ on $U$, and hence
$\theta=0$ on $U$.
Equation~\eqref{eq:pure_trace_K} therefore gives $K=0$ on $U$.
\end{proof}

\begin{Thm}\label{thm:hessian_weak}
Let $(M^n,g,\nabla)$ be a connected Hessian manifold with $n\geq2$. Then $(M,g,\nabla)$ is weakly CHC if and only if it is of CHC $c$ for some constant $c\in\RR$.
\end{Thm}

\begin{proof}
It suffices to prove the ``only if'' part.
Assume that $(M,g,\nabla)$ is weakly CHC.
By Proposition~\ref{prop:CHC_rigidity}~\ref{item:ptCHC-CHC_rigidity}, it suffices to show
that $(M,g,\nabla)$ is pointwise CHC.
Let $a,b$ be the functions in
Proposition~\ref{prop:weak_characterization}.

Taking the trace in $X$ of~\eqref{eq:weak_integrability}
and using $R^\nabla=0$, we obtain
\begin{equation}\label{eq:hessian_contraction}
\begin{aligned}
 2(n-1)b\Kflat(U,Y,Z)
 ={}&\{2da(U)+(n-1)db(U)\}g(Y,Z)\\
 &-da(Y)g(U,Z)-da(Z)g(U,Y).
\end{aligned}
\end{equation}
We interchange $U$ and $Y$ in~\eqref{eq:hessian_contraction}
and subtract the resulting equation from the original one.
By the symmetry of $\Kflat$, this gives
\begin{equation}\label{eq:hessian_coefficient_symmetry}
 \{3da+(n-1)db\}(U)g(Y,Z)
 =\{3da+(n-1)db\}(Y)g(U,Z).
\end{equation}
Taking the $g$-trace in $Y$ and $Z$ of~\eqref{eq:hessian_coefficient_symmetry}
and using $n\geq2$, we obtain
\begin{equation}\label{eq:hessian_coefficient_derivative}
 3da+(n-1)db=0.
\end{equation}
On the open set $\Omega:=\{p\in M\mid b(p)\ne0\}$,
equations~\eqref{eq:hessian_contraction}
and~\eqref{eq:hessian_coefficient_derivative}
give~\eqref{eq:pure_trace} with
$\theta=-da/(2(n-1)b)$.
Lemma~\ref{lem:pure_trace_vanishing} therefore gives $K=0$ on $\Omega$.
Hence $Q=-\nabla K=0$ on $\Omega$.
The uniqueness of the coefficients in Proposition~\ref{prop:weak_characterization} gives $a=b=0$ on $\Omega$.
Since $b\ne0$ on $\Omega$ by definition,
we conclude that $\Omega=\emptyset$.
Thus $b=0$ everywhere, and Proposition~\ref{prop:weak_characterization}
shows that $(M,g,\nabla)$ is pointwise CHC.
\end{proof}

\subsection{Examples of statistical manifolds of CHC}
\label{subsec:non-Hessian_Ex}

In this subsection, we construct non-trivial conjugate symmetric statistical structures of CHC using a warped product.
Here, non-trivial means that $\nabla\ne\nabla^g$ (see Section~\ref{subsec:Stat-str}).

\begin{Prop}\label{prop:warped_CHC}
Let $(N,h)$ be a Riemannian manifold with $\dim N \geq 1$ and
$\mu\ne0$ a constant.
On $M=\RR\times N$ put
\[
 g=dt^2+e^{-2\mu t}h,\qquad \xi=\partial_t,
\]
and define a symmetric $(1,2)$-tensor field $K$ by
\begin{equation}\label{eq:K-def_CS-CHC}
    \begin{aligned}
         K(\xi,\xi)&=2\mu\xi,\qquad K(\xi,X)=\mu X,\\
         K(X,Y)&=\mu g(X,Y)\xi \qquad (X,Y\in TN).
    \end{aligned}
\end{equation}
Then $(g,\nabla=\nabla^g+K)$ is a conjugate symmetric statistical
structure of CHC $4\mu^2$ and constant sectional
$\mathsf K$-curvature $\mu^2$.  
Moreover, for vector fields $X,Y,Z$ on $N$, the curvature of $\nabla$ is given by
\begin{equation}\label{eq:R-nabla_CHC-Ex}
 R^\nabla(\xi,X)=0,\qquad R^\nabla(X,Y)\xi=0,\qquad
 R^\nabla(X,Y)Z=R^h(X,Y)Z.
\end{equation}
In particular, $(M,g,\nabla)$ is Hessian if and only if $h$ is flat.
In the non-Hessian case, $g$ does not have constant sectional curvature.
\end{Prop}

\begin{proof}
By the definition of $K$ in~\eqref{eq:K-def_CS-CHC},
the tensor $g(K(\cdot,\cdot),\cdot)$ is totally symmetric.
Hence $(g,\nabla)$ is a statistical structure.
Let $X,Y,Z$ be vector fields on $N$ lifted independently of $t$.
The Levi-Civita connection formulas for warped products
\cite[Lemma~7.3]{Bishop-ONeill_1969}, together with
$\nabla=\nabla^g+K$, give
\begin{equation}\label{eq:nabla_CHC-Ex}
 \nabla_\xi\xi=2\mu\xi,\qquad
 \nabla_\xi X=\nabla_X\xi=0,\qquad
 \nabla_XY=\nabla^h_XY+2\mu g(X,Y)\xi.
\end{equation}
These formulas, together with $\xi(g(X,Y))=-2\mu g(X,Y)$ and $\nabla^h h = 0$, give
\begin{equation}\label{eq:nabla-K_CHC-Ex}
\begin{aligned}
 (\nabla_\xi K)(\xi,\xi)&=-4\mu^2\xi,\\
 (\nabla_\xi K)(\xi,X)&=-2\mu^2X,\\
 (\nabla_\xi K)(X,Y)&=0,\\
 (\nabla_X K)(\xi,\xi)&=0,\\
 (\nabla_X K)(\xi,Y)&=-2\mu^2g(X,Y)\xi,\\
 (\nabla_X K)(Y,Z)
 &=-2\mu^2\{g(X,Y)Z+g(X,Z)Y\}.
\end{aligned}
\end{equation}
It follows from~\eqref{eq:nabla-K_CHC-Ex} that
\[
 (\nabla_UK)(V,W)
 =-2\mu^2\{g(U,V)W+g(U,W)V\}
\]
for all vector fields $U,V,W$ on $M$.
Since $Q=-\nabla K$, it follows from Corollary~\ref{cor:weak_forces_CS} and
Proposition~\ref{prop:CHC_tensor} that $(M,g,\nabla)$
is conjugate symmetric and of CHC $4\mu^2$.
Moreover, by Proposition~\ref{prop:CHC_curvature_relation},
it has constant sectional $\mathsf K$-curvature $\mu^2$.
The curvature identities in~\eqref{eq:R-nabla_CHC-Ex} follow directly
from~\eqref{eq:nabla_CHC-Ex} and show that $R^\nabla=0$ if and only
if $R^h=0$. 
The last assertion follows from Corollary~\ref{cor:CHC_nonHessian_sectional_curvature}, since $\mu\ne0$.
\end{proof}

\begin{Rem}
When $(N,h)$ is Euclidean space, the structure in
Proposition~\ref{prop:warped_CHC} is statistically isomorphic
to the Hessian manifold of CHC $4\mu^2$ constructed by Shima~\cite[Theorem~7]{Shima_1995}.
\end{Rem}

\section{Left-invariant Hessian structures on \texorpdfstring{$\mathbb R^n$}{Euclidean n-space}}
\label{sec:euclidean}

In this section, we focus on left-invariant Hessian structures on
the abelian Lie group $\RR^n$ and their moduli spaces.
In particular, we compute their principal Hessian curvatures
and show that these invariants determine their statistical
isomorphism classes (see Theorem~\ref{thm:euclidean_spectrum}).

A statistical structure $(g,\nabla)$ on a Lie group $G$ is
\emph{left-invariant} if both $g$ and $\nabla$ are preserved by
every left translation $L_a:G\to G$, $x\mapsto ax$, with $a\in G$
(cf.~\cite{FIK}).

\subsection{Moduli spaces}\label{subsec:left_invariant_moduli}

In this subsection, we fix our terminology for moduli spaces
of left-invariant statistical structures on Lie groups and discuss their relation to principal Hessian curvatures.

Let $G$ be a Lie group and write
\[
 \LStat(G):=\{(g,\nabla) \mid (g,\nabla)\text{ is a left-invariant statistical structure on }G\}.
\]
The automorphism group $\Aut(G)$ acts on $\LStat(G)$ by
\[
 \varphi\cdot(g,\nabla)
 =\bigl((\varphi^{-1})^*g,(\varphi^{-1})^*\nabla\bigr).
\]
We denote the resulting equivalence relation by $\sim_{\HS}$.
The \emph{moduli space of left-invariant statistical structures on $G$}
is defined to be the quotient space
\[
 \mathcal M_{\HS}(G)
 :=\LStat(G)/{\sim_{\HS}}.
\]
Since $(g,\nabla)\sim_{\HS}(g',\nabla')$ implies
$(g,\nabla)\sim_{\STAT}(g',\nabla')$ (see Definition~\ref{def:STAT_equivalence}), there is a natural surjection
\begin{equation}\label{eq:HS_to_STAT}
 \mathcal M_{\HS}(G)\longrightarrow\LStat(G)/{\sim_{\STAT}}.
\end{equation}
We equip $\mathcal M_{\HS}(G)$ with the quotient topology; for the topology on $\LStat(G)$, see \cite[Section~3.2]{KOOT_2025}.
The corresponding moduli spaces of conjugate symmetric and Hessian
structures are denoted by $\mathcal M_{\mathrm{CS}}(G)$ and
$\mathcal M_{\Hess}(G)$, respectively, and carry the relative topology.

\begin{Rem}\label{rem:moduli_no_scaling}
In our previous work~\cite{KOOT_2025}, we defined the moduli space as
\[
 \MLStat(G):=(\RR_{>0}\times\Aut(G))\backslash\LStat(G).
\]
Here we omit the $\RR_{>0}$-action, since principal Hessian curvatures
generally change under metric scaling.
\end{Rem}

Let $\gee$ denote the Lie algebra of $G$, that is, the space of
left-invariant vector fields on $G$.
As is well known, left-invariant tensor fields on $G$
can be naturally identified with tensors on $\gee$
(see~\cite{Kobayashi-Nomizu_II}).
We use the same notation for a left-invariant tensor field and the corresponding tensor on $\gee$.

In particular, for a left-invariant conjugate symmetric
statistical structure $(g,\nabla)$ on $G$, we regard the
Hessian sectional curvature as a function
$q:S^2\gee\setminus\{0\}\to\RR$.
Moreover, the Hessian curvature operator $\hat{Q}^{(g,\nabla)}$
is a self-adjoint operator on $S^2\gee$ with respect to
$\LR_g$.
We call its eigenvalues the \emph{principal Hessian curvatures
of $(g,\nabla)$} and denote their multiset, counted with
multiplicity, by $\PHC(g,\nabla)$.
By left-invariance, $\PHC(g,\nabla)=\PHC_a(g,\nabla)$ for every $a\in G$
(see Definition~\ref{def:Hessian_curvatures} for the notation).

The principal Hessian curvatures are preserved by statistical
isomorphisms (see~Theorem~\ref{thm:HSC-def-on-CS}) and define a continuous map on
$\mathcal M_{\mathrm{CS}}(G)$, as follows.

\begin{Prop}\label{prop:PHC_moduli_map}
Let $G$ be an $n$-dimensional Lie group, and put $N=n(n+1)/2$.
Let $\mathscr S_N:=\RR^N/\mathfrak S_N$ denote the space of unordered
real $N$-tuples, equipped with the quotient topology.
Then the principal Hessian curvatures define a continuous map
\begin{equation}\label{eq:PHC_moduli_map}
 \mathcal P_G:\mathcal M_{\mathrm{CS}}(G)\longrightarrow\mathscr S_N,
 \qquad [(g,\nabla)]\longmapsto\PHC(g,\nabla).
\end{equation}
In particular, statistical structures with distinct values of
$\mathcal P_G$ are not statistically isomorphic.
\end{Prop}

Since the proof of the preceding proposition is straightforward, we omit it here.

\subsection{Preliminaries on \texorpdfstring{$\mathbb R^n$}{Euclidean n-space}}\label{subsec:Rn_Hesse-moduli}

In this subsection, we collect some facts about left-invariant statistical structures on the Lie group $\RR^n$ that will be used later.
Let $x^1,\ldots,x^n$ denote the standard coordinates on $\RR^n$,
and write
\[
 g_E:=\sum_{i=1}^n(dx^i)^2
\]
for the Euclidean metric.

We first discuss some general properties of left-invariant
statistical structures on $\RR^n$.
Let $(g,C)$ be a left-invariant statistical structure on $\RR^n$.
Since $\RR^n$ is abelian, the Koszul formula (cf.~\cite[Chapter~IV, \S~2, Proposition~2.3]{Kobayashi-Nomizu_I}) gives $\nabla^g_XY=0$ for all left-invariant vector fields $X,Y$.
As $C$ is left-invariant, it follows that
\begin{equation}\label{eq:Rn_parallel_tensors}
 \nabla^gC=0.
\end{equation}
In particular, every such structure is conjugate symmetric
(see Section~\ref{subsec:class_structures}).
Thus the principal Hessian
curvatures and Hessian sectional curvature are defined for every left-invariant statistical
structure on $\RR^n$.

We next discuss the moduli space $\mathcal M_{\HS}(\RR^n)$.
Recall that $\sim_{\STAT}$ is the equivalence relation
defined by statistical isomorphism (see Definition~\ref{def:STAT_equivalence}).
Using the same notation for its restriction to $\LStat(\RR^n)$,
we have
\begin{equation}\label{eq:equivalence_Rn}
 \sim_{\HS}\;=\;\sim_{\STAT}
 \qquad\text{on }\LStat(\RR^n).
\end{equation}
Indeed, let $f$ be a statistical isomorphism between two
left-invariant statistical structures $(g,C)$ and $(g',C')$ on $\RR^n$.
Under the standard identifications of tangent spaces with
$\RR^n$, its differential at $0$ defines a linear map
$L\in GL(n,\RR)$ satisfying $L^*g'=g$ and $L^*C'=C$
by left-invariance.
Thus $L\in\Aut(\RR^n)$ is also a statistical isomorphism.

Let $\gee_{\RR^n}$ denote the Lie algebra of $\RR^n$.
Every left-invariant Riemannian metric on $\RR^n$ is equivalent
to $g_E$ under the action of $GL(n,\RR)=\Aut(\RR^n)$.
Since the stabilizer of $g_E$ is $O(n)$, we obtain a bijection
\begin{equation}\label{eq:Rn_statistical_moduli}
 \begin{aligned}
 O(n)\backslash S^3((\gee_{\RR^n})^*)
 &\longrightarrow
 \mathcal M_{\HS}(\RR^n)=\mathcal M_{\mathrm{CS}}(\RR^n),\\
 [C]&\longmapsto[(g_E,C)].
 \end{aligned}
\end{equation}
This bijection is a homeomorphism for the quotient topologies,
where $S^3((\gee_{\RR^n})^*)$ is equipped with its usual
finite-dimensional vector space topology.
For the topological argument, see~\cite[Theorem~A.1]{KOOT_2025}.

We next recall the classification of left-invariant Hessian
structures on $\RR^n$.
Put
\begin{equation}\label{eq:lambda_chamber}
 \Lambda_n:=\{\lambda\in\RR^n\mid
                \lambda_1\geq\cdots\geq\lambda_n\geq0\},
 \qquad
 C_\lambda:=\sum_{i=1}^n\lambda_i(dx^i)^3.
\end{equation}
By~\eqref{eq:equivalence_Rn} and
\cite[Theorems~5.1(2) and~C.1]{KOOT_2025},
we obtain the following.

\begin{Prop}\label{prop:diagonal_normal_form}
Every left-invariant Hessian structure on $\RR^n$ is
statistically isomorphic to $(g_E,C_\lambda)$ for a unique
$\lambda\in\Lambda_n$.
Conversely, $(g_E,C_\lambda)$ is Hessian for every
$\lambda\in\Lambda_n$.
\end{Prop}

\begin{Rem}\label{rem:euclidean_classification_background}
For the relation of this classification to the earlier 
classification of homogeneous Hessian domains,
see~\cite[Remark~5.5]{KOOT_2025}.
\end{Rem}

\subsection{Principal Hessian curvatures of left-invariant Hessian structures}
\label{subsec:euclidean_spectrum}

In this subsection, we give explicit formulas for the Hessian sectional
curvature and principal Hessian curvatures of all left-invariant
Hessian structures on $\RR^n$.

We retain the notation of the previous subsection.
Let $e_i=\partial/\partial x^i$ and put
\begin{equation}\label{eq:euclidean_symmetric_basis}
 F_{ii}=e_i\odot e_i,\qquad
 F_{ij}=\sqrt2\,e_i\odot e_j\quad(i<j).
\end{equation}
These vectors form an orthonormal basis of $S^2\gee_{\RR^n}$
with respect to $\LR_{g_E}$.

\begin{Thm}\label{thm:euclidean_spectrum}
For $\lambda\in\Lambda_n$, the Hessian curvature operator of
$(g_E,C_\lambda)$ satisfies
\begin{equation}\label{eq:euclidean_A}
 \hat{Q}(F_{ii}) = \frac{\lambda_i^2}{4}F_{ii},\qquad
 \hat{Q}(F_{ij}) = 0\quad(i<j).
\end{equation}
Its principal Hessian curvatures, counted with multiplicities, are
\begin{equation}\label{eq:euclidean_spectrum}
 \frac{\lambda_1^2}{4},\ldots,\frac{\lambda_n^2}{4},
 \quad\underbrace{0,\ldots,0}_{n(n-1)/2}.
\end{equation}
For $0\ne\nu=\sum_{i\leq j}\nu_{ij}F_{ij}$,
\begin{equation}\label{eq:euclidean_q}
 q (\nu)=
 \frac{\sum_{i=1}^n\lambda_i^2\nu_{ii}^2}
 {4\sum_{i\leq j}\nu_{ij}^2}.
\end{equation}
In particular, two left-invariant Hessian structures on $\RR^n$ are
statistically isomorphic if and only if their principal Hessian
curvatures agree, with multiplicities.
\end{Thm}

\begin{proof}
Let $K$ be the difference tensor of $(g_E,C_\lambda)$.
Then
\[
 K(e_i,e_i)=-\frac{\lambda_i}{2}e_i,\qquad
 K(e_i,e_j)=0\quad(i\ne j).
\]
Thus $[K,K]=0$, and~\eqref{eq:Rn_parallel_tensors} implies
$\nabla^{g_E}K=0$.
Lemma~\ref{lem:gram_curvature}, applied with $\rho=0$, gives
\[
 \calQ(\nu,\mu)=g_E(\calK\nu,\calK\mu).
\]
Hence $\hat{Q} = \calK^* \circ \calK$. Since
\[
 \calK(F_{ii})=-\frac{\lambda_i}{2}e_i,\qquad
 \calK(F_{ij})=0\quad(i<j),
\]
we obtain~\eqref{eq:euclidean_A} and~\eqref{eq:euclidean_spectrum}.
Equation~\eqref{eq:euclidean_q} then follows from~\eqref{eq:Rayleigh}.
For the last assertion, the ``only if'' part follows from
Theorem~\ref{thm:HSC-def-on-CS}.
The converse follows from Proposition~\ref{prop:diagonal_normal_form}
and~\eqref{eq:euclidean_spectrum}.
\end{proof}

Combining Theorems~\ref{thm:euclidean_spectrum} and~\ref{thm:hessian_weak}, we obtain the following corollary.

\begin{Cor}\label{cor:euclidean_hessian_weak}
If $n\geq2$, a left-invariant Hessian structure on $\RR^n$ is weakly
CHC if and only if it is trivial.
\end{Cor}

\begin{Rem}\label{rem:homogeneous_dimension_one}
For $n=1$ and $\lambda>0$, the left-invariant Hessian structure
$(g_E,\lambda(dx)^3)$ is non-trivial and, by Theorem~\ref{thm:euclidean_spectrum},
is of CHC $\lambda^2/4$.
Thus the dimension assumption in Corollary~\ref{cor:euclidean_hessian_weak} is essential.
\end{Rem}

\section{Weakly CHC statistical structures with parallel cubic forms}\label{sec:parallel}

\subsection{Classification}\label{subsec:WCHC_Classification}
In this subsection, we present a classification of weakly CHC
statistical structures in dimensions at least $2$ whose cubic forms
are parallel with respect to the Levi-Civita connection $\nabla^g$ (Theorem~\ref{thm:parallel_classification}).
Since $\nabla^gC=0$ implies that $\nabla^gC$ is totally symmetric, these structures are conjugate symmetric (see Section~\ref{subsec:class_structures}).

The non-trivial models in the classification below arise from
classical affine spheres, which we now describe.
We follow~\cite{Nomizu-Sasaki} for the terminology of affine
differential geometry and refer to it for details.
For $c>0$, set
\begin{equation}\label{eq:sigma_model}
 \Sigma_c=\{(x,y,z)\in(\RR_{>0})^3\mid xyz=c\}.
\end{equation}
The inclusion $\Sigma_c\hookrightarrow\RR^3$ is a locally
strongly convex immersion.
We choose the sign of its affine normal field so that the
induced affine metric is positive definite.
Let $g_c$ and $\nabla_c$ denote the resulting affine metric
and induced affine connection (Blaschke connection),
respectively.
The pair $(g_c,\nabla_c)$ is a statistical structure on $\Sigma_c$.
The surface $\Sigma_c$ is a hyperbolic affine sphere
(see Dumas--Wolf~\cite[Section~4.4]{Dumas-Wolf_2015} and the original work of \c{T}i\c{t}eica~\cite[Section~9]{Tzitzeica_1908}).

\begin{Thm}\label{thm:parallel_classification}
Let $(M^n,g,\nabla)$ be a connected statistical manifold with $n\geq2$,
and suppose that its cubic form $C = \nabla g$ satisfies $\nabla^gC=0$.
Then $(M,g,\nabla)$ is weakly CHC if and only if one of the following holds:
\begin{enumerate}[label=\textup{(\arabic*)}]
 \item\label{item:parallel_trivial} $\nabla=\nabla^g$.
 \item\label{item:parallel_sigma} $n=2$, and $(M,g,\nabla)$ is locally
 statistically isomorphic to $(\Sigma_c,g_c,\nabla_c)$ for some $c>0$.
\end{enumerate}
In the latter case~\ref{item:parallel_sigma}, $c$ is uniquely determined, and if $M$ is simply connected and $g$ is complete, then $(M,g,\nabla)$ is
globally statistically isomorphic to $(\Sigma_c,g_c,\nabla_c)$.
\end{Thm}

The proof is given in Appendix~\ref{subsec:parallel_classification_proof}.

\begin{Cor}
Let $(M^n,g,\nabla)$ be a connected, simply connected, and weakly CHC
statistical manifold of dimension $n\geq2$ with complete metric $g$.
If $\nabla^gC=0$ and $C\not\equiv0$, then $n=2$ and
$(M,g,\nabla)$ is statistically isomorphic to
$(\Sigma_c,g_c,\nabla_c)$ for a unique $c>0$.
\end{Cor}

The statistical manifold $(\Sigma_c,g_c,\nabla_c)$ can be
described in terms of a left-invariant statistical structure
on $\RR^2$.
We now give this description.
For $\alpha\geq0$, let $K^{(\alpha)}$ be the left-invariant $(1,2)$-tensor field on $\RR^2$ defined by
\begin{equation}\label{eq:parallel_surface_normal}
\begin{aligned}
 K^{(\alpha)}(e_1,e_1)&=\alpha e_1,\\
 K^{(\alpha)}(e_1,e_2)&=K^{(\alpha)}(e_2,e_1)=-\alpha e_2,\\
 K^{(\alpha)}(e_2,e_2)&=-\alpha e_1
\end{aligned}
\end{equation}
in the standard frame $e_i=\partial/\partial x^i$ $(i=1,2)$.
Put
\begin{equation}\label{eq:euclidean_parallel_model}
 \nabla^{(\alpha)}=\nabla^{g_E}+K^{(\alpha)},\qquad
 C^{(\alpha)}=\alpha\bigl(-2(dx^1)^3+6\,dx^1(dx^2)^2\bigr).
\end{equation}
Then $(g_E,\nabla^{(\alpha)})$ is a left-invariant statistical structure with cubic form $C^{(\alpha)}$.

\begin{Prop}\label{prop:parallel_models}
For every $\alpha\geq0$, the statistical structure
$(g_E,\nabla^{(\alpha)})$ is trace-free and $\nabla^{g_E} C^{(\alpha)} = 0$.
Moreover,
\begin{equation}\label{eq:parallel_model_invariants}
 \lVert C^{(\alpha)}\rVert_{g_E}=4\alpha,\qquad
 R^{\nabla^{(\alpha)}}=-2\alpha^2R_0,
\end{equation}
and its Hessian curvature operator is
\begin{equation}\label{eq:parallel_surface_spectrum}
 \hat{Q} = 4\alpha^2\pr_{\tr} - 2\alpha^2\pr_0.
\end{equation}
In particular, it is weakly CHC, and it is of CHC if and only if $\alpha=0$.
For $\alpha>0$, $(\RR^2, g_E, \nabla^{(\alpha)})$ is globally statistically isomorphic to $(\Sigma_{c(\alpha)},g_{c(\alpha)},\nabla_{c(\alpha)})$, where
\begin{equation}\label{eq:sigma_parameter}
 c(\alpha)=\frac{1}{4\cdot3^{3/2}\alpha^4}.
\end{equation}
Consequently, $(\Sigma_c,g_c)$ is complete, simply connected, and flat,
and $\nabla^{g_c}C_c=0$, where $C_c=\nabla_cg_c$.
\end{Prop}

The proof of Proposition~\ref{prop:parallel_models} is given in
Appendix~\ref{subsec:parallel_model_proof}.

We conclude this subsection with a remark on the assumption
$\nabla^gC=0$ in Theorem~\ref{thm:parallel_classification}.
For a statistical manifold $(M,g,C)$, define the \emph{divergence} of $C$ by
\begin{equation}\label{eq:cubic_divergence}
 (\operatorname{div}_gC)(X,Y)
 =\sum_{i=1}^n(\nabla^g_{e_i}C)(e_i,X,Y),
\end{equation}
where $\{ e_1,\ldots,e_n \}$ is any local $g$-orthonormal frame.
If conjugate symmetry is assumed, the assumption $\nabla^gC=0$
in Theorem~\ref{thm:parallel_classification} can be replaced
by any of the equivalent conditions in
Proposition~\ref{prop:weak_parallel_equivalences} below.

\begin{Prop}\label{prop:weak_parallel_equivalences}
Let $(M^n,g,\nabla)$ be a weakly CHC conjugate symmetric
statistical manifold with $n\geq2$.
The following conditions are equivalent:
\begin{enumerate}[label=\textup{(\roman*)}]
\item\label{item:weak_parallel_C}
$\nabla^gC=0$;
\item\label{item:weak_divergence_zero}
$\operatorname{div}_gC=0$;
\item\label{item:weak_divergence_puretrace}
$\operatorname{div}_gC=\phi g$ for some function $\phi:M\to\RR$;
\item\label{item:weak_tracefree}
$\tr_gK=0$.
\end{enumerate}
In particular, the function $\phi$
in~\ref{item:weak_divergence_puretrace} must vanish.
\end{Prop}

The proof is given in Appendix~\ref{subsec:parallel_divergence_proof}.
No connectedness assumption is needed in this proposition.
For $n\geq3$, the four conditions are also equivalent to $C=0$,
by Theorem~\ref{thm:parallel_classification} applied to each
connected component.

\subsection{Moduli space of left-invariant weakly CHC statistical structures on \texorpdfstring{$\mathbb R^n$}{Euclidean n-space}}

In this subsection, we apply the classification result from the
previous subsection to left-invariant statistical structures on
$\RR^n$ and describe their moduli space.

\begin{Cor}\label{cor:Rn_weak_classification}
Every left-invariant weakly CHC statistical structure on $\RR^2$ is
statistically isomorphic to $(g_E,\nabla^{(\alpha)})$ for a unique
$\alpha\in[0,\infty)$.
For $n\geq3$, every left-invariant weakly CHC statistical structure
on $\RR^n$ is statistically isomorphic to $(g_E,\nabla^{g_E})$.
\end{Cor}

\begin{proof}
Since every left-invariant Riemannian metric on $\RR^n$ is complete,
the existence of the asserted statistical isomorphisms follows from
\eqref{eq:Rn_parallel_tensors},
Theorem~\ref{thm:parallel_classification}, and
Proposition~\ref{prop:parallel_models}.
Uniqueness of $\alpha$ follows from the first identity
in~\eqref{eq:parallel_model_invariants} and the invariance
of $\lVert C\rVert_g$ under statistical isomorphisms.
\end{proof}

Let $\mathcal M_{\mathrm{wCHC}}(\RR^n)$ denote the moduli space
of left-invariant weakly CHC statistical structures on $\RR^n$,
equipped with the relative topology from
$\mathcal M_{\HS}(\RR^n)$, as in Section~\ref{subsec:left_invariant_moduli}.

\begin{Cor}\label{cor:Rn_weak_moduli}
The map
\[
 [0,\infty)\longrightarrow\mathcal M_{\mathrm{wCHC}}(\RR^2),\qquad
 \alpha\longmapsto[(g_E,C^{(\alpha)})]
\]
is a homeomorphism. For $n\geq3$,
$\mathcal M_{\mathrm{wCHC}}(\RR^n)$ consists of a single point.
\end{Cor}

\begin{proof}
By \eqref{eq:equivalence_Rn} and
Corollary~\ref{cor:Rn_weak_classification}, the map is bijective.
It is continuous by the quotient topology, since $C^{(\alpha)}$
depends continuously on $\alpha$.
Its inverse is $[(g,C)]\mapsto\frac14\lVert C\rVert_g$,
by the first identity in~\eqref{eq:parallel_model_invariants}.
This inverse is continuous because $(g,C)\mapsto\lVert C\rVert_g$
is a continuous $\Aut(\RR^2)$-invariant function on $\LStat(\RR^2)$.
The assertion for $n\geq3$ follows from
\eqref{eq:equivalence_Rn} and Corollary~\ref{cor:Rn_weak_classification}.
\end{proof}

\section*{Acknowledgements}
The author would like to express his sincere gratitude to Ryu Ueno
for his valuable suggestions and encouragement, and to Takayuki Okuda
for his insightful suggestions.
The author also thanks Hitoshi Furuhata, Naoto Satoh, Hiroshi Tamaru, and Koichi Tojo for their valuable comments and suggestions.

\section*{Use of AI}
The author used ChatGPT to revise the English and
to help locate some of the references.
ChatGPT also suggested using results from affine differential
geometry to identify the local model in the non-trivial case
of Theorem~\ref{thm:parallel_classification}.
The author verified the details of this approach
and briefly discusses it in
Remark~\ref{rem:affine-resuls_for_classification}
at the end of Appendix~\ref{subsec:parallel_classification_proof}.
The author takes full responsibility for the content of this paper.

\appendix

\section{Proofs for Section~\ref{sec:parallel}}\label{app:parallel}

\subsection{Proof of Proposition~\ref{prop:parallel_models}}\label{subsec:parallel_model_proof}

In this subsection, we give a proof of Proposition~\ref{prop:parallel_models}.

\begin{proof}[Proof of Proposition~\ref{prop:parallel_models}]
The equality $\nabla^{g_E}C^{(\alpha)}=0$ follows
from~\eqref{eq:Rn_parallel_tensors}.
In particular, the structure is conjugate symmetric
by~\eqref{eq:CS_Codazzi}.
Direct calculation using
\eqref{eq:parallel_surface_normal} and
\eqref{eq:euclidean_parallel_model} gives
\[
 \tr_{g_E}K^{(\alpha)}=0,\qquad
 \lVert C^{(\alpha)}\rVert_{g_E}=4\alpha,\qquad
 [K^{(\alpha)},K^{(\alpha)}]=-2\alpha^2R_0.
\]
Since $g_E$ is flat, \eqref{eq:conjugate_curvature} gives
the curvature identity in~\eqref{eq:parallel_model_invariants}.
Applying~\eqref{eq:gram_operator} in
Lemma~\ref{lem:gram_curvature} with $\rho=-2\alpha^2$,
we obtain~\eqref{eq:parallel_surface_spectrum}.
By Propositions~\ref{prop:weak_characterization}
and~\ref{prop:CHC_tensor},
the structure is weakly CHC for every $\alpha\geq0$, and it is of CHC if and only if $\alpha=0$.

We now establish the statistical isomorphism.
We follow~\cite{Nomizu-Sasaki} for the terminology
of affine differential geometry.
For $\alpha>0$, put $c=c(\alpha)$ as
in~\eqref{eq:sigma_parameter}, and define
\begin{equation}\label{eq:sigma_explicit_map}
 \Phi(u,v)=\Phi_\alpha(u,v):=c^{1/3}
 \bigl(e^{-\alpha u+\sqrt3\alpha v},
       e^{-\alpha u-\sqrt3\alpha v},e^{2\alpha u}\bigr),
\end{equation}
where $(u,v)=(x^1,x^2)$ are the standard coordinates on $\RR^2$.
This is a diffeomorphism from $\RR^2$ onto $\Sigma_c$.
We show that $\Phi$ is a statistical isomorphism from
$(\RR^2,g_E,\nabla^{(\alpha)})$ to $(\Sigma_c,g_c,\nabla_c)$.

Regard $\Phi$ as an immersion into $\RR^3$ and define
a vector field $\xi$ along $\Phi$ by
\[
 \xi_{(u,v)}:=2\alpha^2\Phi(u,v),
\]
using the standard identification
$T_{\Phi(u,v)}\RR^3\simeq\RR^3$.
We show that $\xi$ is the affine normal field of $\Phi$.
Direct calculation gives
\begin{equation}\label{eq:eq-volume}
 \det(\Phi_*\partial_u,\Phi_*\partial_v,\xi)
 =12\sqrt3\,c\alpha^4=1.
\end{equation}
Thus, at every $(u,v)\in\RR^2$,
\[
 T_{\Phi(u,v)}\RR^3
 =\Phi_*(T_{(u,v)}\RR^2)\oplus\RR\xi_{(u,v)},
\]
so $\xi$ is transversal to the immersion $\Phi$.

Let $\overline D$ denote the standard flat affine connection
on $\RR^3$.
A direct calculation gives
\begin{align}
 \overline D_X(\Phi_*Y)
 &=\Phi_*(\nabla^{(\alpha)}_X Y)+g_E(X,Y)\xi,
 \label{eq:Phi_gauss}\\
 \overline D_X\xi
 &=2\alpha^2\Phi_*X
 \label{eq:Phi_weingarten}
\end{align}
for all vector fields $X,Y$ on $\RR^2$.
By~\eqref{eq:Phi_gauss}, the affine fundamental form
and affine connection induced by $\xi$ are
$g_E$ and $\nabla^{(\alpha)}$, respectively.
In particular, the affine fundamental form is positive definite.
Equation~\eqref{eq:Phi_weingarten} shows that
$\overline D_X\xi$ is tangent to $\Sigma_c$.
Moreover, \eqref{eq:eq-volume} gives
\[
 \det(\Phi_*X,\Phi_*Y,\xi)=\omega_{g_E}(X,Y),
\]
where $\omega_{g_E}=du\wedge dv$ is the Riemannian volume
form of $g_E$ with respect to the standard orientation on $\RR^2$.
These conditions characterize the affine normal field;
see~\cite[Chapter~II]{Nomizu-Sasaki}.
Thus $\xi$ is the affine normal field of $\Phi$.

By uniqueness of the affine normal field inducing a positive definite
affine metric, the affine normal field of the inclusion
$\Sigma_c\hookrightarrow\RR^3$ takes the value $\xi_{(u,v)}$
at $\Phi(u,v)$. Equation~\eqref{eq:Phi_gauss} therefore gives
\[
 \Phi^*g_c=g_E,\qquad
 \Phi^*\nabla_c=\nabla^{(\alpha)}.
\]
Consequently, $\Phi$ is a global statistical isomorphism from
$(\RR^2,g_E,\nabla^{(\alpha)})$ to $(\Sigma_c,g_c,\nabla_c)$.
The remaining assertions follow immediately from this isomorphism.
\end{proof}

\subsection{Proof of Theorem~\ref{thm:parallel_classification}}
\label{subsec:parallel_classification_proof}

In this subsection, we discuss the proof of
Theorem~\ref{thm:parallel_classification}.

By the statistical isomorphism in
Proposition~\ref{prop:parallel_models} and the bijectivity of
$\alpha\mapsto c(\alpha)$ in~\eqref{eq:sigma_parameter}
on $(0,\infty)$, it suffices to prove the following proposition.

\begin{Prop}\label{prop:classification_Rn-left-inv}
Let $(M^n,g,\nabla)$ be a connected statistical manifold
with $n\geq2$ and $\nabla^gC=0$.
Then $(M,g,\nabla)$ is weakly CHC if and only if
one of the following holds:
\begin{enumerate}[label=\textup{(\arabic*)}]
 \item\label{item:parallel_trivial_inAp}
 $\nabla=\nabla^g$.
 \item\label{item:Rn-isom_inAp}
 $n=2$, and $(M,g,\nabla)$ is locally statistically
 isomorphic to $(\RR^2,g_E,\nabla^{(\alpha)})$
 for some $\alpha>0$.
\end{enumerate}
In case~\ref{item:Rn-isom_inAp}, the constant $\alpha$ is uniquely
determined, and if $M$ is simply connected and $g$ is complete,
then there exists a global statistical isomorphism
from $(M,g,\nabla)$ to $(\RR^2,g_E,\nabla^{(\alpha)})$.
\end{Prop}

We first prove the ``if'' part of
Proposition~\ref{prop:classification_Rn-left-inv}.
If $\nabla=\nabla^g$, then $K=0$ and hence $Q=-\nabla K=0$. 
Thus the structure is of CHC $0$ and hence weakly CHC.
The statistical manifolds $(\RR^2,g_E,\nabla^{(\alpha)})$ are weakly CHC
by Proposition~\ref{prop:parallel_models}.
By Theorem~\ref{thm:HSC-def-on-CS} and
Definition~\ref{def:weak_CHC}, the weakly CHC condition
is preserved under the relation of being locally
statistically isomorphic.
Thus case~\ref{item:Rn-isom_inAp} also implies that
$(M,g,\nabla)$ is weakly CHC.

For the ``only if'' part, we prepare the following two lemmas.

\begin{Lem}\label{lem:weak_parallel_dimension_tracefree}
Let $(M^n,g,\nabla)$ be a weakly CHC statistical
manifold with $n\geq2$ and $\nabla^gC=0$. 
If $K\not\equiv0$, then $n=2$ and $\tr_gK=0$.
\end{Lem}

\begin{Lem}\label{lem:weak_parallel_surface_geometry}
Let $(M^2,g,\nabla)$ be a connected trace-free weakly CHC
statistical manifold with $\nabla^gC=0$ and $K\not\equiv0$.
Then $g$ is flat, and there is a constant $\alpha>0$
such that, for every $p\in M$, there exists a
$g_p$-orthonormal basis $\{e_1,e_2\}$ of $T_pM$ satisfying
\begin{equation}\label{eq:weak_parallel_surface_coefficients}
 \begin{aligned}
 K_p(e_1,e_1)&=\alpha e_1,\\
 K_p(e_1,e_2)&=-\alpha e_2,\\
 K_p(e_2,e_2)&=-\alpha e_1.
 \end{aligned}
\end{equation}
Moreover,
\begin{equation}\label{eq:weak_parallel_alpha_norm}
 \alpha=\frac14\lVert C\rVert_g.
\end{equation}
\end{Lem}

We also use the following result of Satoh.

\begin{Prop}[{\cite[Corollary~3.5]{Satoh-2019}}]
\label{prop:Satoh_local_isomorphism}
Let $(M,g,\nabla)$ and $(M',g',\nabla')$ be statistical
manifolds of the same dimension, with difference tensors
$K$ and $K'$, respectively.
Suppose that
\[
 \nabla^gR^g=0,\qquad \nabla^gK=0,\qquad
 \nabla^{g'}R^{g'}=0,\qquad \nabla^{g'}K'=0.
\]
Let $p\in M$ and $p'\in M'$, and let
$T:T_pM\to T_{p'}M'$ be a linear isometry satisfying
\[
 \begin{aligned}
 T\bigl(R^g_p(u,v)w\bigr)
 &=R^{g'}_{p'}(Tu,Tv)Tw,\\
 T\bigl(K_p(u,v)\bigr)
 &=K'_{p'}(Tu,Tv)
 \end{aligned}
\]
for all $u,v,w\in T_pM$.
Then the map
\[
 F=\exp_{p'}^{g'}\circ T\circ(\exp_p^g)^{-1}
\]
defines a statistical isomorphism between sufficiently small neighborhoods of $p$ and $p'$.
\end{Prop}

Assuming the two lemmas above, we now prove the ``only if''
part of Proposition~\ref{prop:classification_Rn-left-inv},
together with the uniqueness and global assertions.
The proofs of the lemmas are given afterwards.

\begin{proof}[Completion of the proof of
Proposition~\ref{prop:classification_Rn-left-inv}]
Suppose that $(M,g,\nabla)$ is weakly CHC.
If $K\equiv0$, then $\nabla=\nabla^g$, giving
case~\ref{item:parallel_trivial_inAp}.
In what follows, we assume that $K\not\equiv0$.

By Lemma~\ref{lem:weak_parallel_dimension_tracefree},
$n=2$ and the structure is trace-free.
Lemma~\ref{lem:weak_parallel_surface_geometry} shows that
$g$ is flat and provides a constant $\alpha>0$.
Fix an arbitrary point $p\in M$.
By this lemma, we may choose a $g_p$-orthonormal basis
$\{e'_1,e'_2\}$ of $T_pM$ such that
\[
 K_p(e'_1,e'_1)=\alpha e'_1,\qquad
 K_p(e'_1,e'_2)=-\alpha e'_2,\qquad
 K_p(e'_2,e'_2)=-\alpha e'_1.
\]
Let $e_i=\partial/\partial x^i$ $(i=1,2)$ be the standard
frame on $\RR^2$, and define $T:T_pM\to T_0\RR^2$ by
\[
 T(e'_i)=(e_i)_0\qquad(i=1,2).
\]
By construction, $T$ is a linear isometry, and
\eqref{eq:parallel_surface_normal} gives
\[
 T\bigl(K_p(u,v)\bigr)
 =K^{(\alpha)}_0(Tu,Tv)
 \qquad(u,v\in T_pM).
\]

The assumption $\nabla^gC=0$ gives $\nabla^gK=0$
by~\eqref{eq:K-C}.
By Proposition~\ref{prop:parallel_models}, we also have
$\nabla^{g_E}K^{(\alpha)}=0$.
Both $g$ and $g_E$ are flat.
Proposition~\ref{prop:Satoh_local_isomorphism} therefore
implies that, on a sufficiently small connected normal
neighborhood $U$ of $p$, the map
\begin{equation}\label{eq:parallel_local_isomorphism}
 F:=\exp_0^{g_E}\circ T\circ(\exp_p^g)^{-1}
\end{equation}
is a statistical isomorphism onto an open neighborhood
of $0\in\RR^2$, with $F(p)=0$ and $dF_p=T$.
Since $p$ was arbitrary and $(g_E,\nabla^{(\alpha)})$ is
left-invariant, $(M,g,\nabla)$ is locally statistically
isomorphic to $(\RR^2,g_E,\nabla^{(\alpha)})$.

We next prove the uniqueness of $\alpha$.
Suppose that $(M,g,\nabla)$ is also locally statistically
isomorphic to $(\RR^2,g_E,\nabla^{(\alpha')})$
for some $\alpha'>0$.
By Proposition~\ref{prop:parallel_models} and invariance
of the norm of the cubic form under statistical isomorphisms,
$4\alpha=4\alpha'$.
Hence $\alpha=\alpha'$.

Finally, suppose that $M$ is simply connected and $g$ is complete.
By the Cartan--Hadamard theorem (cf.~\cite[Chapter~VIII, Theorem~8.1]{Kobayashi-Nomizu_II})
and the flatness of $g$, \eqref{eq:parallel_local_isomorphism} defines a global isometry $F:M\to\RR^2$.
Since $F^*K^{(\alpha)}$ and $K$ are $\nabla^g$-parallel
and agree at $p$, they coincide on the connected manifold $M$.
Thus $F^*\nabla^{(\alpha)}=\nabla$, and $F$ is the required
global statistical isomorphism.
\end{proof}

We first prove Lemma~\ref{lem:weak_parallel_dimension_tracefree}.
For the proof, we prepare the following lemma.

\begin{Lem}\label{lem:weak_divergence_trace}
Let $(M^n,g,\nabla)$ be a weakly CHC conjugate symmetric
statistical manifold with $n\geq2$.
Let $a,b$ be the functions in
Proposition~\ref{prop:weak_characterization}, and put
\begin{equation}\label{eq:weak_divergence_notation}
 E:=\tr_gK,\qquad
 \tau:=g(E,\cdot),\qquad
 P:=(n+1)a+b.
\end{equation}
Then
\begin{align}
 \operatorname{div}_gC
 &=-2\nabla^g\tau,
 \label{eq:weak_divergence_tau}\\
 (\nabla^g_X\tau)(Y)
 &=\Kflat(E,X,Y)-P\,g(X,Y)
 \label{eq:weak_divergence_trace}
\end{align}
for all vector fields $X,Y$ on $M$.
\end{Lem}

\begin{proof}
By~\eqref{eq:trace_free_equivalences}, we have
$\tau(X)=\tr K_X$ and $\tr_g C(X,\cdot,\cdot)=-2\tau(X)$.
Conjugate symmetry implies that $\nabla^gC$ is totally
symmetric.
Thus, for a local orthonormal frame $\{ e_1,\ldots,e_n \}$,
\eqref{eq:cubic_divergence} gives
\[
 \begin{aligned}
 (\operatorname{div}_gC)(X,Y)
 &=\sum_{i=1}^n(\nabla^g_{e_i}C)(e_i,X,Y)\\
 &=\sum_{i=1}^n(\nabla^g_XC)(Y,e_i,e_i)\\
 &=-2(\nabla^g_X\tau)(Y).
 \end{aligned}
\]
This proves~\eqref{eq:weak_divergence_tau}.
Taking the trace of the endomorphism
$Z\mapsto(\nabla^g_XK)(Y,Z)$
in~\eqref{eq:weak_CHC_nablaK}, we obtain
\[
 \begin{aligned}
 (\nabla^g_X\tau)(Y)
 &=\tau(K(X,Y))-\{na+(a+b)\}g(X,Y)\\
 &=\Kflat(E,X,Y)-P\,g(X,Y),
 \end{aligned}
\]
which proves~\eqref{eq:weak_divergence_trace}.
\end{proof}

We now prove Lemma~\ref{lem:weak_parallel_dimension_tracefree}.

\begin{proof}[Proof of Lemma~\ref{lem:weak_parallel_dimension_tracefree}]
Let $a,b$ be the functions in
Proposition~\ref{prop:weak_characterization},
and let $E,\tau,P$ be as in~\eqref{eq:weak_divergence_notation}.
Put $\rho:=(a-b)/2$.
By~\eqref{eq:K-C}, the assumption $\nabla^gC=0$
gives $\nabla^gK=0$, and hence $\nabla^g\tau=0$.
Lemma~\ref{lem:weak_K_curvature} gives
$[K,K]=\rho R_0$.
Since $K$ and $R_0$ are $\nabla^g$-parallel,
the identity $[K,K]=\rho R_0$ gives $d\rho=0$.
Equation~\eqref{eq:weak_CHC_drho} then gives $b\tau=0$.
Since $\nabla^g\tau=0$,
\eqref{eq:weak_divergence_trace}
and~\eqref{eq:trace_free_equivalences} give
\begin{equation}\label{eq:parallel_trace_identity}
 K_E=P\id,\qquad \lVert E\rVert_g^2 = \tr K_E=nP.
\end{equation}
Suppose that $E_p\ne0$ at some point $p\in M$.
Then $b(p)=0$ by $b\tau=0$.
Choose a unit vector $X\in T_pM$ perpendicular to $E_p$.
Using $[K,K]=\rho R_0$ and~\eqref{eq:parallel_trace_identity},
we obtain
\[
 0=[K_{E_p},K_X]X=\rho(p)E_p.
\]
Thus $\rho(p)=0$, and hence $a(p)=0$.
Consequently, $P(p)=0$, contradicting
$|E_p|^2=nP(p)$.
Therefore $\tr_gK=E=0$, and
\eqref{eq:parallel_trace_identity} also gives $P=0$ on $M$.

It remains to show that $n=2$.
Suppose, for a contradiction, that $n\geq3$.
Since $K\not\equiv0$, choose $p\in M$ such that $K_p\ne0$.
Equation~\eqref{eq:weak_CHC_bK} gives $b(p)=0$, since $n\geq3$ and $K_p\ne0$.
The identity $P(p)=0$ then gives $a(p)=0$, and
\eqref{eq:weak_tensor} gives $Q_p=0$.
Equations~\eqref{eq:Q_expansion} and~\eqref{eq:B_K},
together with $\nabla^gK=0$, give
\[
 0=Q_p(X,X,X,X)=|K_p(X,X)|^2
 \qquad(X\in T_pM).
\]
Polarization yields $K_p=0$, a contradiction.
Thus $n=2$.
\end{proof}

We next prove Lemma~\ref{lem:weak_parallel_surface_geometry}.
For the proof, we use the following result of Opozda.

\begin{Lem}[{\cite[Theorem~4.4]{Opozda_2016}}]
\label{lem:Opozda_tracefree_constant_K}
Let $(M,g,\nabla)$ be a trace-free conjugate symmetric
statistical manifold. 
If the structure has constant sectional $\mathsf K$-curvature, then $\nabla^gK=0$.
Moreover, if $K\not\equiv0$, then $g$ is flat.
\end{Lem}

\begin{proof}[Proof of Lemma~\ref{lem:weak_parallel_surface_geometry}]
By conjugate symmetry and the weakly CHC condition, Lemma~\ref{lem:weak_K_curvature} gives $[K,K]=\rho R_0$ for a function $\rho$.
Since $K$ and $R_0$ are $\nabla^g$-parallel and $M$ is connected,
$\rho$ is constant.
Since $K\not\equiv0$, Lemma~\ref{lem:Opozda_tracefree_constant_K} implies that $g$ is flat.
Since $K$ is $\nabla^g$-parallel and $M$ is connected,
the assumption $K\not\equiv0$ implies that $K$ is nowhere zero.
For each $p\in M$, the trace-free case of
\cite[Lemma~3.6]{Opozda_2016}, applied with $n=2$,
gives an orthonormal basis of $T_pM$ in which $K_p$
has the form~\eqref{eq:weak_parallel_surface_coefficients}
with a parameter $\alpha_p>0$.
The relation $C=-2\Kflat$ gives
\[
 \alpha_p=\frac14|C_p|.
\]
Since $C$ is $\nabla^g$-parallel and $M$ is connected,
$|C_p|$ is independent of $p$.
This completes the proof.
\end{proof}

\begin{Rem}\label{rem:affine-resuls_for_classification}
In the ``only if'' part of Theorem~\ref{thm:parallel_classification},
the local identification in case~\ref{item:parallel_sigma}
can also be established using results from affine differential geometry
after Lemmas~\ref{lem:weak_parallel_dimension_tracefree}
and~\ref{lem:weak_parallel_surface_geometry} have been proved.
In the non-trivial case, these lemmas give $n=2$, $\tr_gK=0$,
and flatness of $g$.
Equations~\eqref{eq:conjugate_curvature}
and~\eqref{eq:weak_parallel_surface_coefficients} then give
$R^\nabla=-2\alpha^2R_0$ for a constant $\alpha>0$.
Thus, in particular, the statistical structure is trace-free,
and its dual connection $\nabla^*$ is projectively flat.
With some further argument using the realization theorem of
Dillen--Nomizu--Vranken~\cite[Theorem~3, p.~229]{Vranken+_1990}
and the classification of proper affine spheres with positive
definite flat affine metric by
Magid--Ryan~\cite[Theorem~2, p.~283]{Magid-Ryan_1990},
one can show that $(M,g,\nabla)$ and
$(\Sigma_c,g_c,\nabla_c)$ are locally statistically isomorphic
for some $c>0$.
\end{Rem}

\subsection{Proof of Proposition~\ref{prop:weak_parallel_equivalences}}
\label{subsec:parallel_divergence_proof}

In this subsection, we prove
Proposition~\ref{prop:weak_parallel_equivalences}.

\begin{proof}[Proof of Proposition~\ref{prop:weak_parallel_equivalences}]
The implications
\ref{item:weak_parallel_C}$\Rightarrow$\ref{item:weak_divergence_zero}
$\Rightarrow$\ref{item:weak_divergence_puretrace}
are immediate.
Let $a,b$ be the functions in
Proposition~\ref{prop:weak_characterization},
and let $E,\tau,P$ be as in~\eqref{eq:weak_divergence_notation}.
Put $\rho:=(a-b)/2$.

We prove~\ref{item:weak_divergence_puretrace}
$\Rightarrow$\ref{item:weak_tracefree}.
Assume~\ref{item:weak_divergence_puretrace}, and put
$s:=-\phi/2$ and $\beta:=s+P$.
Lemma~\ref{lem:weak_divergence_trace} gives
\begin{equation}\label{eq:nabla-g_sX_K_E}
 \nabla^g_XE=sX,\qquad K_E=\beta\id.
\end{equation}
Taking the trace of the second identity
in~\eqref{eq:nabla-g_sX_K_E}
and using~\eqref{eq:trace_free_equivalences}, we obtain
\begin{equation}\label{eq:weak_equiv_beta_trace}
 n\beta=\tr K_E=\tau(E)=\lVert E\rVert_g^2.
\end{equation}
We show that $U:=\{p\in M:E_p\ne0\}$ is empty.
For $p\in U$, choose a unit vector $X\in T_pM$
perpendicular to $E_p$.
By~\eqref{eq:weak_rho}
and the second identity in~\eqref{eq:nabla-g_sX_K_E},
we have
\begin{equation}\label{eq:weak_equiv_commutator}
 0=[K_{E_p},K_X]X=\rho(p)E_p.
\end{equation}
Since $E_p\ne0$ and $p\in U$ was arbitrary,
$\rho=0$ on $U$.
As $U$ is open, $d\rho=0$ on $U$, and
\eqref{eq:weak_CHC_drho} gives $b\tau=0$.
Since $\tau$ is nowhere zero on $U$, it follows that $b=0$.
The definition of $\rho$ then gives $a=0$, and hence $P=0$.
Therefore, $\beta=s$, and~\eqref{eq:weak_equiv_beta_trace} gives
\begin{equation}\label{eq:weak_equiv_beta_positive}
 \beta=s=\frac{ \lVert E\rVert_g^2 }{n}>0
 \qquad\text{on }U.
\end{equation}
Differentiating the second identity
in~\eqref{eq:nabla-g_sX_K_E} with respect to $\nabla^g_Z$
and using the first identity, we obtain
\begin{equation}\label{eq:nabla-g_K_d-beta}
 (\nabla^g_ZK)(E,Y)+sK(Z,Y)=d\beta(Z)Y.
\end{equation}
By~\eqref{eq:CS_Codazzi} and the symmetry of $K$,
the left-hand side of~\eqref{eq:nabla-g_K_d-beta}
is symmetric in $Y,Z$.
Therefore, we have
\begin{equation}\label{eq:weak_equiv_dbeta_symmetry}
 d\beta(Z)Y=d\beta(Y)Z.
\end{equation}
Taking the trace in $Y$
in~\eqref{eq:weak_equiv_dbeta_symmetry}, we obtain
$(n-1)d\beta(Z)=0$.
Since $n\geq2$, we have $d\beta=0$ on $U$.
Differentiating~\eqref{eq:weak_equiv_beta_trace}
and using the first identity in~\eqref{eq:nabla-g_sX_K_E}
now yields
\begin{equation}\label{eq:weak_equiv_contradiction}
 0=n\,d\beta(Z)
 =2g(\nabla^g_ZE,E)
 =2s\tau(Z)
\end{equation}
for every vector field $Z$ on $U$.
Taking $Z=E$ in~\eqref{eq:weak_equiv_contradiction}
gives $0=2s \lVert E\rVert_g^2$, contradicting
\eqref{eq:weak_equiv_beta_positive} if $U$ is nonempty.
Hence $U$ is empty, so $E=0$.
This proves~\ref{item:weak_divergence_puretrace}
$\Rightarrow$\ref{item:weak_tracefree}.

Finally, we prove~\ref{item:weak_tracefree}
$\Rightarrow$\ref{item:weak_parallel_C}.
Assume~\ref{item:weak_tracefree}.
Then $\tau=0$, and~\eqref{eq:weak_CHC_drho} gives
$d\rho=0$.
Thus, by~\eqref{eq:weak_rho}, the structure has constant
sectional $\mathsf K$-curvature on each connected component.
Since the structure is trace-free and conjugate symmetric,
Lemma~\ref{lem:Opozda_tracefree_constant_K}, applied to each connected component, gives $\nabla^gK=0$.
Equation~\eqref{eq:K-C} then gives $\nabla^gC=0$,
proving~\ref{item:weak_tracefree}
$\Rightarrow$\ref{item:weak_parallel_C}.
The assertion $\phi=0$ follows from
\ref{item:weak_divergence_puretrace}$\Rightarrow$\ref{item:weak_divergence_zero}.
\end{proof}

\providecommand{\bysame}{\leavevmode\hbox to3em{\hrulefill}\thinspace}
\providecommand{\href}[2]{#2}
\ifdefined\nolinkurl
  \providecommand{\amsplainid}[1]{{\urlstyle{same}%
    \setbox0=\hbox{\nolinkurl{#1}}%
    \ifdim\wd0>.75\linewidth \nolinkurl{#1}\else\box0\fi}}
\else
  \providecommand{\amsplainid}[1]{{\ttfamily\escapechar=-1\relax\detokenize{#1}}}
\fi


\begin{thebibliography}{10}

\bibitem{Amari-Nagaoka_2000}
Shun-ichi Amari and Hiroshi Nagaoka, \emph{Methods of information geometry}, Translations of Mathematical Monographs, vol. 191, American Mathematical Society and Oxford University Press, 2000.
\href{https://doi.org/10.1090/mmono/191}{doi:~\amsplainid{10.1090/mmono/191}}.

\bibitem{Bishop-ONeill_1969}
Richard~L. Bishop and Barrett O'Neill, \emph{Manifolds of negative curvature}, Trans. Amer. Math. Soc. \textbf{145} (1969), 1--49.
\href{https://doi.org/10.1090/S0002-9947-1969-0251664-4}{doi:~\amsplainid{10.1090/S0002-9947-1969-0251664-4}}.

\bibitem{Calabi_1972}
Eugenio Calabi, \emph{Complete affine hyperspheres. {I}}, in \emph{Symposia {M}athematica, {V}ol. {X}}, Academic Press, 1972, pp.~19--38.
\href{https://www.ams.org/mathscinet-getitem?mr=365607}{MR~\amsplainid{365607}}.

\bibitem{Chen-Mihai-Mihai_2019}
Bang-Yen Chen, Adela Mihai, and Ion Mihai, \emph{A {C}hen first inequality for statistical submanifolds in {H}essian manifolds of constant {H}essian curvature}, Results Math. \textbf{74} (2019), no.~4, Art. no.~165.
\href{https://doi.org/10.1007/s00025-019-1091-y}{doi:~\amsplainid{10.1007/s00025-019-1091-y}}.

\bibitem{Vranken+_1990}
Franki Dillen, Katsumi Nomizu, and Luc Vrancken, \emph{Conjugate connections and {R}adon's theorem in affine differential geometry}, Monatsh. Math. \textbf{109} (1990), no.~3, 221--235.
\href{https://doi.org/10.1007/BF01297762}{doi:~\amsplainid{10.1007/BF01297762}}.

\bibitem{Dumas-Wolf_2015}
David Dumas and Michael Wolf, \emph{Polynomial cubic differentials and convex polygons in the projective plane}, Geom. Funct. Anal. \textbf{25} (2015), no.~6, 1734--1798.
\href{https://doi.org/10.1007/s00039-015-0344-5}{doi:~\amsplainid{10.1007/s00039-015-0344-5}}.

\bibitem{Furuhata_2009}
Hitoshi Furuhata, \emph{Hypersurfaces in statistical manifolds}, Differential Geom. Appl. \textbf{27} (2009), no.~3, 420--429.
\href{https://doi.org/10.1016/j.difgeo.2008.10.019}{doi:~\amsplainid{10.1016/j.difgeo.2008.10.019}}.

\bibitem{Furuhata_2011}
\bysame, \emph{Statistical hypersurfaces in the space of {H}essian curvature zero}, Differential Geom. Appl. \textbf{29} (2011), S86--S90.
\href{https://doi.org/10.1016/j.difgeo.2011.04.012}{doi:~\amsplainid{10.1016/j.difgeo.2011.04.012}}.

\bibitem{Furuhata-Hasegawa-Satoh_2022}
Hitoshi Furuhata, Izumi Hasegawa, and Naoto Satoh, \emph{Chen invariants and statistical submanifolds}, Commun. Korean Math. Soc. \textbf{37} (2022), no.~3, 851--864.
\href{https://ckms.kms.or.kr/journal/view.html?doi=10.4134/CKMS.c210185}{doi:~\amsplainid{10.4134/CKMS.c210185}}.

\bibitem{FIK}
Hitoshi Furuhata, Jun-ichi Inoguchi, and Shimpei Kobayashi, \emph{A characterization of the alpha-connections on the statistical manifold of normal distributions}, Inf. Geom. \textbf{4} (2021), no.~1, 177--188.
\href{https://doi.org/10.1007/s41884-020-00037-z}{doi:~\amsplainid{10.1007/s41884-020-00037-z}}.

\bibitem{FK_2013}
Hitoshi Furuhata and Takashi Kurose, \emph{Hessian manifolds of nonpositive constant {H}essian sectional curvature}, Tohoku Math. J. (2) \textbf{65} (2013), no.~1, 31--42.
\href{https://doi.org/10.2748/tmj/1365452623}{doi:~\amsplainid{10.2748/tmj/1365452623}}.

\bibitem{Inoguchi-Ohno_2024}
Jun-ichi Inoguchi and Yu~Ohno, \emph{Homogeneous statistical manifolds}, Inf. Geom. \textbf{8} (2025), no.~2, 285--341.
\href{https://doi.org/10.1007/s41884-025-00172-5}{doi:~\amsplainid{10.1007/s41884-025-00172-5}}.

\bibitem{IKO_2026}
Hideyuki Ishi, Hikozo Kobayashi, and Takayuki Okuda, \emph{Invariant statistical connections on the multivariate centered Gaussian model and their moduli spaces}, 2026.
\href{https://arxiv.org/abs/2605.30300}{arXiv:\amsplainid{2605.30300}}.

\bibitem{Khan-Zhang_2020}
Gabriel Khan and Jun Zhang, \emph{The {K}\"ahler geometry of certain optimal transport problems}, Pure Appl. Anal. \textbf{2} (2020), no.~2, 397--426.
\href{https://doi.org/10.2140/paa.2020.2.397}{doi:~\amsplainid{10.2140/paa.2020.2.397}}.

\bibitem{Khan-Zhang-Zheng_2025}
Gabriel Khan, Jun Zhang, and Fangyang Zheng, \emph{Positively curved {K}\"ahler metrics on tube domains and their applications to optimal transport}, Math. Res. Lett. \textbf{32} (2025), no.~2, 457--491.
\href{https://doi.org/10.4310/mrl.250721011640}{doi:~\amsplainid{10.4310/mrl.250721011640}}.

\bibitem{Khan-Zheng_2026}
Gabriel Khan and Fangyang Zheng, \emph{K\"ahler-{R}icci flow preserves negative anti-bisectional curvature}, Trans. Amer. Math. Soc. \textbf{379} (2026), no.~8, 5319--5351.
\href{https://doi.org/10.1090/tran/9584}{doi:~\amsplainid{10.1090/tran/9584}}.

\bibitem{Kirillov_2008}
Alexander Kirillov, Jr., \emph{An Introduction to {Lie} Groups and {Lie} Algebras}, Cambridge Studies in Advanced Mathematics, vol. 113, Cambridge University Press, 2008.
\href{https://doi.org/10.1017/CBO9780511755156}{doi:~\amsplainid{10.1017/CBO9780511755156}}.

\bibitem{KOOT_2025}
Hikozo Kobayashi, Yu~Ohno, Takayuki Okuda, and Hiroshi Tamaru, \emph{The moduli spaces of left-invariant statistical structures on {L}ie groups}, Inf. Geom. \textbf{9} (2026), no.~1, 235--278.
\href{https://doi.org/10.1007/s41884-026-00200-y}{doi:~\amsplainid{10.1007/s41884-026-00200-y}}.

\bibitem{Kobayashi-Okuda_2025}
Hikozo Kobayashi and Takayuki Okuda, \emph{On invariant conjugate symmetric statistical structures on the space of zero-mean multivariate normal distributions}, in \emph{Geometric science of information. {P}art {I}}, Lecture Notes in Comput. Sci., vol. 16033, Springer, 2026, pp.~65--72.
\href{https://doi.org/10.1007/978-3-032-03918-7_7}{doi:~\amsplainid{10.1007/978-3-032-03918-7_7}}.

\bibitem{KO-CC}
Shimpei Kobayashi and Yu~Ohno, \emph{On a constant curvature statistical manifold}, Inf. Geom. \textbf{5} (2022), no.~1, 31--46.
\href{https://doi.org/10.1007/s41884-022-00065-x}{doi:~\amsplainid{10.1007/s41884-022-00065-x}}.

\bibitem{Kobayashi-Ohno_2025}
\bysame, \emph{A characterization of the alpha-connections on the statistical manifold of multivariate normal distributions}, Osaka J. Math. \textbf{62} (2025), no.~2, 329--349.
\href{https://doi.org/10.18910/101132}{doi:~\amsplainid{10.18910/101132}}.

\bibitem{Kobayashi-Nomizu_I}
Shoshichi Kobayashi and Katsumi Nomizu, \emph{Foundations of differential geometry. {V}ol {I}}, Interscience Publishers, 1963.
\href{https://www.ams.org/mathscinet-getitem?mr=152974}{MR~\amsplainid{152974}}.

\bibitem{Kobayashi-Nomizu_II}
\bysame, \emph{Foundations of differential geometry. {V}ol. {II}}, Interscience Tracts in Pure and Applied Mathematics, no.~15, Interscience Publishers, 1969.
\href{https://www.ams.org/mathscinet-getitem?mr=238225}{MR~\amsplainid{238225}}.

\bibitem{Kurose-1990}
Takashi Kurose, \emph{Dual connections and affine geometry}, Math. Z. \textbf{203} (1990), no.~1, 115--121.
\href{https://doi.org/10.1007/BF02570725}{doi:~\amsplainid{10.1007/BF02570725}}.

\bibitem{Kurose_1994}
\bysame, \emph{On the divergences of {$1$}-conformally flat statistical manifolds}, Tohoku Math. J. (2) \textbf{46} (1994), no.~3, 427--433.
\href{https://doi.org/10.2748/tmj/1178225722}{doi:~\amsplainid{10.2748/tmj/1178225722}}.

\bibitem{L-SM}
Steffen~L. Lauritzen, \emph{Statistical manifolds}, in \emph{Differential geometry in statistical inference}, Institute of Mathematical Statistics Lecture Notes--Monograph Series, vol.~10, Institute of Mathematical Statistics, 1987, pp.~163--216.
\href{https://doi.org/10.1214/lnms/1215467061}{doi:~\amsplainid{10.1214/lnms/1215467061}}.

\bibitem{Li-Penn_1988}
An~Min Li and Gabi Penn, \emph{Uniqueness theorems in affine differential geometry. {II}}, Results Math. \textbf{13} (1988), no.~3-4, 308--317.
\href{https://doi.org/10.1007/BF03323247}{doi:~\amsplainid{10.1007/BF03323247}}.

\bibitem{Loftin_2010}
John Loftin, \emph{Survey on affine spheres}, in \emph{Handbook of geometric analysis, {N}o. 2}, Adv. Lect. Math. (ALM), vol.~13, Int. Press, 2010, pp.~161--191.
\href{https://www.ams.org/mathscinet-getitem?mr=2743442}{MR~\amsplainid{2743442}}.

\bibitem{Loftin-McIntosh_2016}
John Loftin and Ian McIntosh, \emph{Cubic differentials in the differential geometry of surfaces}, in \emph{Handbook of {T}eichm\"uller theory. {V}ol. {VI}}, IRMA Lect. Math. Theor. Phys., vol.~27, Eur. Math. Soc., 2016, pp.~231--274.
\href{https://doi.org/10.4171/161-1/9}{doi:~\amsplainid{10.4171/161-1/9}}.

\bibitem{Magid-Nomizu_1989}
Martin~A. Magid and Katsumi Nomizu, \emph{On affine surfaces whose cubic forms are parallel relative to the affine metric}, Proc. Japan Acad. Ser. A Math. Sci. \textbf{65} (1989), no.~7, 215--218.
\href{https://doi.org/10.3792/pjaa.65.215}{doi:~\amsplainid{10.3792/pjaa.65.215}}.

\bibitem{Magid-Ryan_1990}
Martin~A. Magid and Patrick~J. Ryan, \emph{Flat affine spheres in {${\bf R}^3$}}, Geom. Dedicata \textbf{33} (1990), no.~3, 277--288.
\href{https://doi.org/10.1007/BF00181334}{doi:~\amsplainid{10.1007/BF00181334}}.

\bibitem{Matsuzoe-Takeuchi-Amari_2006}
Hiroshi Matsuzoe, Jun-ichi Takeuchi, and Shun-ichi Amari, \emph{Equiaffine structures on statistical manifolds and {B}ayesian statistics}, Differential Geom. Appl. \textbf{24} (2006), no.~6, 567--578.
\href{https://doi.org/10.1016/j.difgeo.2006.02.003}{doi:~\amsplainid{10.1016/j.difgeo.2006.02.003}}.

\bibitem{Mihai-Mihai_2020}
Adela Mihai and Ion Mihai, \emph{The {$\delta(2, 2)$}-invariant on statistical submanifolds in {H}essian manifolds of constant {H}essian curvature}, Entropy \textbf{22} (2020), no.~2, Art. no.~164.
\href{https://doi.org/10.3390/e22020164}{doi:~\amsplainid{10.3390/e22020164}}.

\bibitem{Mihai_2020}
Ion Mihai, \emph{Statistical manifolds and their submanifolds. {R}esults on {C}hen-like invariants}, in \emph{Geometry of submanifolds}, Contemp. Math., vol. 756, Amer. Math. Soc., 2020, pp.~163--172.
\href{https://doi.org/10.1090/conm/756/15206}{doi:~\amsplainid{10.1090/conm/756/15206}}.

\bibitem{Min-Ri-Kwak-An_2015}
Cholrim Min, Wonhak Ri, Kumhyok Kwak, and Dokjun An, \emph{Equiaffine structure and conjugate {R}icci-symmetry of a statistical manifold}, Differential Geom. Appl. \textbf{41} (2015), 39--47.
\href{https://doi.org/10.1016/j.difgeo.2015.04.005}{doi:~\amsplainid{10.1016/j.difgeo.2015.04.005}}.

\bibitem{Molitor_2022}
Mathieu Molitor, \emph{One-dimensional exponential families with constant {H}essian sectional curvature}, Inf. Geom. \textbf{5} (2022), no.~2, 511--530.
\href{https://doi.org/10.1007/s41884-022-00074-w}{doi:~\amsplainid{10.1007/s41884-022-00074-w}}.

\bibitem{Figueiredo-Molitor_2026}
Danuzia Nascimento~Figueir\^edo and Mathieu Molitor, \emph{Toric dually flat manifolds and complex space forms}, Differential Geom. Appl. \textbf{102} (2026), Art. no.~102321.
\href{https://doi.org/10.1016/j.difgeo.2025.102321}{doi:~\amsplainid{10.1016/j.difgeo.2025.102321}}.

\bibitem{Noguchi_1992}
Mitsunori Noguchi, \emph{Geometry of statistical manifolds}, Differential Geom. Appl. \textbf{2} (1992), no.~3, 197--222.
\href{https://doi.org/10.1016/0926-2245(92)90011-B}{doi:~\amsplainid{10.1016/0926-2245(92)90011-B}}.

\bibitem{Nomizu-Sasaki}
Katsumi Nomizu and Takeshi Sasaki, \emph{Affine differential geometry: Geometry of affine immersions}, Cambridge Tracts in Mathematics, vol. 111, Cambridge University Press, 1994.
\href{https://www.ams.org/mathscinet-getitem?mr=1311248}{MR~\amsplainid{1311248}}.

\bibitem{Opozda_2015}
Barbara Opozda, \emph{Bochner's technique for statistical structures}, Ann. Global Anal. Geom. \textbf{48} (2015), no.~4, 357--395.
\href{https://doi.org/10.1007/s10455-015-9475-z}{doi:~\amsplainid{10.1007/s10455-015-9475-z}}.

\bibitem{Opozda_2016}
\bysame, \emph{A sectional curvature for statistical structures}, Linear Algebra Appl. \textbf{497} (2016), 134--161.
\href{https://doi.org/10.1016/j.laa.2016.02.021}{doi:~\amsplainid{10.1016/j.laa.2016.02.021}}.

\bibitem{Opozda_2019}
\bysame, \emph{Curvature bounded conjugate symmetric statistical structures with complete metric}, Ann. Global Anal. Geom. \textbf{55} (2019), no.~4, 687--702.
\href{https://doi.org/10.1007/s10455-019-09647-y}{doi:~\amsplainid{10.1007/s10455-019-09647-y}}.

\bibitem{Opozda_2021}
\bysame, \emph{Completeness in affine and statistical geometry}, Ann. Global Anal. Geom. \textbf{59} (2021), no.~3, 367--383.
\href{https://doi.org/10.1007/s10455-021-09752-x}{doi:~\amsplainid{10.1007/s10455-021-09752-x}}.

\bibitem{Opozda_2022}
\bysame, \emph{Some inequalities and applications of {S}imons' type formulas in {R}iemannian, affine, and statistical geometry}, J. Geom. Anal. \textbf{32} (2022), no.~4, Art. no.~108.
\href{https://doi.org/10.1007/s12220-021-00828-5}{doi:~\amsplainid{10.1007/s12220-021-00828-5}}.

\bibitem{Opozda_2025}
\bysame, \emph{The geometry of the tangent and sphere bundles over statistical manifolds}, Results Math. \textbf{80} (2025), no.~2, Art. no.~60.
\href{https://doi.org/10.1007/s00025-025-02360-2}{doi:~\amsplainid{10.1007/s00025-025-02360-2}}.

\bibitem{Sakamoto_2026}
Hakobi Sakamoto, \emph{Hessian manifolds of positive constant Hessian sectional curvature}, 2026.
\href{https://arxiv.org/abs/2610.01737}{arXiv:\amsplainid{2610.01737}}.

\bibitem{Satoh_2007}
Hiroyasu Satoh, \emph{Almost {H}ermitian structures on tangent bundles}, in \emph{Proceedings of the {E}leventh {I}nternational {W}orkshop on {D}ifferential {G}eometry}, Kyungpook Nat. Univ., 2007, pp.~105--118.
\href{https://www.ams.org/mathscinet-getitem?mr=2333235}{MR~\amsplainid{2333235}}.

\bibitem{Satoh-2019}
Naoto Satoh, \emph{Local existence of statistical diffeomorphisms}, JP Journal of Geometry and Topology \textbf{22} (2019), no.~1, 73--95.
\href{https://doi.org/10.17654/GT022010073}{doi:~\amsplainid{10.17654/GT022010073}}.

\bibitem{Shima_1976}
Hirohiko Shima, \emph{On certain locally flat homogeneous manifolds of solvable {L}ie groups}, Osaka J. Math. \textbf{13} (1976), no.~2, 213--229.
\href{https://doi.org/10.18910/12000}{doi:~\amsplainid{10.18910/12000}}.

\bibitem{Shima_1980}
\bysame, \emph{Homogeneous {H}essian manifolds}, Ann. Inst. Fourier (Grenoble) \textbf{30} (1980), no.~3, 91--128.
\href{https://doi.org/10.5802/aif.794}{doi:~\amsplainid{10.5802/aif.794}}.

\bibitem{Shima_1995}
\bysame, \emph{Hessian manifolds of constant {H}essian sectional curvature}, J. Math. Soc. Japan \textbf{47} (1995), no.~4, 735--753.
\href{https://doi.org/10.2969/jmsj/04740735}{doi:~\amsplainid{10.2969/jmsj/04740735}}.

\bibitem{Shima}
\bysame, \emph{The geometry of {H}essian structures}, World Scientific, 2007.
\href{https://doi.org/10.1142/9789812707536}{doi:~\amsplainid{10.1142/9789812707536}}.

\bibitem{Siddiqui-Ahmad-Ozel_2019}
Aliya~Naaz Siddiqui, Kamran Ahmad, and Cenap \"Ozel, \emph{Inequalities for statistical submanifolds in {H}essian manifolds of constant {H}essian curvature}, in \emph{Geometric science of information}, Lecture Notes in Comput. Sci., vol. 11712, Springer, 2019, pp.~389--397.
\href{https://doi.org/10.1007/978-3-030-26980-7_40}{doi:~\amsplainid{10.1007/978-3-030-26980-7_40}}.

\bibitem{Takeuchi-Amari_2005}
Jun-ichi Takeuchi and Shun-ichi Amari, \emph{{$\alpha$}-parallel prior and its properties}, IEEE Trans. Inform. Theory \textbf{51} (2005), no.~3, 1011--1023.
\href{https://doi.org/10.1109/TIT.2004.842703}{doi:~\amsplainid{10.1109/TIT.2004.842703}}.

\bibitem{Tzitzeica_1908}
Georges Tzitz{\'e}ica, \emph{Sur une nouvelle classe de surfaces}, Rend. Circ. Mat. Palermo \textbf{25} (1908), 180--187.
\href{https://doi.org/10.1007/BF03029121}{doi:~\amsplainid{10.1007/BF03029121}}.

\bibitem{Vaisman_2014}
Izu Vaisman, \emph{Hessian geometry on {L}agrange spaces}, Int. J. Math. Math. Sci. (2014), Art. no.~793473.
\href{https://doi.org/10.1155/2014/793473}{doi:~\amsplainid{10.1155/2014/793473}}.

\bibitem{Vrancken-Li-Simon_1991}
Luc Vrancken, An~Min Li, and Udo Simon, \emph{Affine spheres with constant affine sectional curvature}, Math. Z. \textbf{206} (1991), no.~4, 651--658.
\href{https://doi.org/10.1007/BF02571370}{doi:~\amsplainid{10.1007/BF02571370}}.

\bibitem{Weyl_1939}
Hermann Weyl, \emph{The {C}lassical {G}roups. {T}heir {I}nvariants and {R}epresentations}, Princeton University Press, 1939.
\href{https://www.ams.org/mathscinet-getitem?mr=255}{MR~\amsplainid{255}}.

\bibitem{Zhang-Khan_2020}
Jun Zhang and Gabriel Khan, \emph{Statistical mirror symmetry}, Differential Geom. Appl. \textbf{73} (2020), Art. no.~101678.
\href{https://doi.org/10.1016/j.difgeo.2020.101678}{doi:~\amsplainid{10.1016/j.difgeo.2020.101678}}.

\end{thebibliography}
\end{document}